\documentclass[smallcondensed]{svjour3}     
\smartqed  
\usepackage{graphicx}
\usepackage{epstopdf}
\usepackage{natbib}
\usepackage{amsmath}
\usepackage{amssymb}
\usepackage{psfig}
\journalname{Methodology and Computing in Applied Probability}
\begin{document}
\bibliographystyle{plainnat}

\title{Tauberian and Abelian theorems for long-range dependent random fields \vspace{-7cm}
}

\titlerunning{Tauberian and
Abelian theorems}        

\author{Nikolai Leonenko        \and
        Andriy Olenko 
}


\institute{N. Leonenko  \at
              School of Mathematics, Cardiff University,\
Senghennydd Road, Cardiff CF24 4AG, United Kingdom \\
              Tel.: +44-29-20875521\\
              Fax: +44-29-20874199\\
              \email{LeonenkoN@cardiff.ac.uk}                  
           \and
           A. Olenko \at
              Department of Mathematics and Statistics,
La Trobe University,  Victoria, 3086, Australia\\
 Tel.: +61-3-94792609 \\
              Fax: +61-3-94792466\\
              \email{a.olenko@latrobe.edu.au}
}

\date{\it *Will appear in Methodology and Computing in Applied Probability.  The final publication is available at link.springer.com. DOI: 10.1007/s11009-012-9276-9.}

\maketitle

\begin{abstract}
This paper surveys Abelian and Tauberian theorems for long-range dependent random fields.  We describe a framework for asymptotic behaviour of covariance functions or variances of averaged functionals of random fields at infinity and spectral densities at zero. The use of the theorems and their limitations are demonstrated through applications to some new and less-known examples of covariance functions of long-range dependent random fields.
\end{abstract}

\section{Introduction}

In this paper we study asymptotic properties of spectral and covariance functions of random fields. We investigate the connection between (a) the behaviour of covariance functions or variances of averaged functionals of random fields at infinity and (b) the behaviour of spectral distribution functions at zero. In the terminology of the integral transforms the statement "(b) implies (a)" is a theorem of Abelian type and the statement "(a) implies (b)" is a Tauberian theorem.

Abelian and Tauberian theorems are not only of their own interest but also have numerous applications in asymptotic problems of probability theory and statistics (see e.g. surveys in \citealt{bingol, leo2, yak, bin5}).

We concentrate our studies on homogeneous random fields with long-range dependence (LRD), which appear in various applications in signal processing, physics, telecommunications, hydrology, etc. The reader can find more details in  \citealt{leo2,gne,lav2,ber} and the references thereinafter. The volume by \citealt{dou1} contains outstanding surveys of the field. In particular, that volume discusses different definitions of LRD of stationary processes in terms of the autocorrelation function (the integral of the correlation function diverges) or the spectrum (the spectral density has a singularity at zero).

We would like to emphasize that, perhaps to the astonishment of
some readers, multidimensional Abelian and Tauberian theory is more complicated than its one-dimensional counterpart for characteristic functions, and the most results were obtained only in the last two decades.

Abelian and Tauberian theorems for random fields with singular spectrum are of great importance in the theory of limit theorems for functionals of dependent random variables.  They are a key starting point in proofs of non-Gaussian asymptotic behaviour of LRD summands.  Several applications to non-central limit theorems for functionals of LRD random fields can be found in  \citealt{leo1,leo2,dou1}.

Abelian and Tauberian theorems also have been exploited in various statistical applications. Namely, the LRD phenomenon exhibits itself through a slow decay of covariances at infinity or non-integrable covariance functions. However, the majority of statistical methods for LRD models, for example, tests for the presence of a LRD component in real data,  are based on properties of periodograms/estimators of spectral densities. These methods have been greatly developed recently, and rely on a local  specification of the spectral density in a neighborhood of zero.

The aims of the paper are:
 \begin{itemize}
   \item to collect and adjust known Abelian and Tauberian theorems to the case of LRD random fields. Similar asymptotic results were obtained in different ares: functional analysis, theory of integral transforms, probability theory, etc. Some of these results appeared in less-known journals or dissertations. Consequently, this literature is not readily available, and thus it is not properly quoted in scholarly publications. Therefore, accurate translations of these results to the terminology of the probability theory and comparisons among themselves are important problems.
   \item  to show the variety of generalizations for different functionals and functional classes. Each Abelian and Tauberian theorem in the paper is a nontrivial extension to new cases. We present motivating examples, which demonstrate the use of the theorems and their limitations. The examples also show sharpness of results, i.e. that the theorems do not hold under wider assumptions.
   \item  to produce some new examples of covariance and spectral distribution functions of LRD random fields, which have closed-form representations. These results can be used in applied spatial statistical modeling.
 \end{itemize}
The outline of the paper is as follows. In Section 2 we recall the basic definitions and formulae of the spectral theory of random fields. Section 3 introduces two classes of regularly varying functions. Bingham's classical  asymptotic results are presented in Section 4. In Section 5 we discuss Abelian and Tauberian  theorems for covariance functions of LRD random fields. Some extensions for the variances of averaged functionals of random fields are given in Section 6.
The case of $O$-regularly varying asymptotic behaviour is discussed in Section 7. Section 8 contains extensions of isotropic results to "radial directional" homogeneous random fields.

All plots and the verification of computations in the examples were performed using \textit{Maple 15}. Constants in the examples  were chosen to apparently display properties of covariance and spectral functions in plots.

\section{Homogeneous and isotropic random fields}

\begin{definition} A real-valued random field $\xi (x),$ $x\in \mathbb{R}^{n},$ satisfying $\mathrm{E}\xi ^{2}(x)<\infty $ is called homogeneous (in the wide sense) if
its first moment $m(x):=\mathrm{E}\xi (x)$ and covariance
function $B(x,y):=\mathrm{Cov} (\xi (x),\xi (y))$ are invariant with respect
to the group $G=(\mathbb{R}^{n},+)$ of shifts in $\mathbb{R}^{n\text{
}},$ i.e. $m(x)=m(x+z)$, $B(x,y)=B(x+z,\ y+z)$ for any $x,y,z \in
\mathbb{R}^{n}.$
\end{definition}
In other words, $m(x)=\emph{const}$ and the
covariance function $B(x,y)$ depends only on the difference $x-y,$ i.e., $%
B(x,y)=B(x-y).$

The following characterization is due to \citealt{boch}:

\begin{proposition} A function $B(x)$
is the covariance function of a mean-square continuous homogeneous random field $\xi (x)$
if and only if there exists a finite measure $F(\cdot)$ on $(\mathbb{R}%
^{n},\mathcal{B}(\mathbb{R}^{n}))$ such that
\begin{equation}
B(x)=\int_{\mathbb{R}^{n}}e^{i\left\langle x,u \right\rangle
}F(du )=\int_{\mathbb{R}^{n}}\cos \left\langle x,u \right\rangle
F(du )  \label{cov1}
\end{equation}%
with $F(\mathbb{R}^{n})=B(0)<\infty,$ where $\left\langle \cdot,\cdot\right\rangle$ is the inner product in $\mathbb{R}^{n}.$
\end{proposition}

The representation (\ref{cov1}) is
the spectral decomposition of the covariance function, and $F\left(\cdot
\right) $ is the spectral measure of the field $\xi (x)$.

Any spectral measure admits a Lebesgue decomposition into absolutely
continuous, discrete and singular components. If the last two components are absent then the spectral measure is absolutely continuous and
\begin{equation*}
F(\Delta )=\int_{\Delta}\ f(u )\ du ,\quad \mbox{for any}\
\Delta \in \mathcal{B}(\mathbb{R}^{n}).
\end{equation*}%
The function $f(u )$, which is integrable over $\mathbb{R}^{n}$, is
the spectral density function of the random field. If the spectral density
exists, then the spectral decomposition~(\ref{cov1}) can be written as
\begin{equation}\label{Bx}
B(x)=\int_{\mathbb{R}^{n}}e^{i\left\langle x,u \right\rangle
}f(u )\ du .
\end{equation}%
If $B(x)\in L_{1}(\mathbb{R}^{n})$, then clearly $f(u )$
exists, and
$$
f(u )=\frac{1}{(2\pi )^{n}}\int_{\mathbb{R}^{n}}e^{-i\left\langle
x,u \right\rangle }B(x)dx,\quad u \in \mathbb{R}^{n}.
$$

A rotation of the Euclidean space $\mathbb{R}^{n}$ is defined as a transformation $\varrho$ of this space that does not change its orientation and
preserves the distance of the points from the origin: $\left\Vert
\varrho x\right\Vert =\left\Vert x\right\Vert $. The rotations of $\mathbb{R}^{n}$
generate the group $SO(n)$.

\begin{definition} A real-valued random field $\xi (x)$, $x\in \mathbb{R}^{n}$ satisfying $%
\mathrm{E}\xi ^{2}(x)<\infty $ is called homogeneous and isotropic (in the
wide sense) if its mean $m(x)$ is
constant and the covariance function depends only on the Euclidean distance $%
\rho \left( x,y\right) =\left\Vert x-y\right\Vert $ between $x$ and $y,$
i.e., $B(x,y)=B(\left\Vert x-y\right\Vert )$.
\end{definition}
In other words, its
expectation $m(x)$ and covariance function $B(x,y)$ are
invariant with respect to shifts, rotations and reflections in $\mathbb{R}%
^{n}$:%
\begin{equation*}
m(x)=m(x+z),\quad B(x,y)=B(x+z,y+z),\quad  m(x)=m(\varrho x),\quad \mbox{and}
\end{equation*}%
\begin{equation*}
B(x,y)=B(\varrho x,\varrho y),\quad \mbox{for any}\ x,y,z\in \mathbb{R}^{n},\quad \varrho \in
SO(n).
\end{equation*}%
The spectral measure $F\left(\cdot
\right) $ of a homogeneous isotropic random field is
invariant with respect to the group $SO(n)$, i.e., $F(\Delta )=F(\varrho \Delta )$
for every $\varrho \in SO(n)$, $\Delta \in \mathcal{B}(\mathbb{R}^{n}).$

The following covariance function is originally due to \citealt{scho}  and
\citealt{yag}.

\begin{proposition} A function $B(r)$, $r:=\left\Vert x\right\Vert,$
is the covariance function of a mean-square continuous homogeneous isotropic
random field $\xi (x)$, $x\in \mathbb{R}^{n},$ if and only if there exists a
finite measure $G\left(\cdot
\right)$ on $(\mathbb{R}_{+},\mathcal{B}(\mathbb{R}_{+}))$ such
that
\begin{equation}
B(r)=\int_{0}^{\infty }Y_{n}(\lambda r)\ G(d\lambda),  \label{cov4}
\end{equation}%
where $Y_{n}(z)$ is the spherical Bessel function defined by%
\begin{equation*}
Y_{1}(z):=\cos z,
\end{equation*}%
\begin{equation*}
Y_{n}(z):=2^{(n-2)/2}\Gamma \left(\frac{n}{2}\right)\ J_{(n-2)/2}(z)\ z^{(2-n)/2},\quad
z\geq 0,n\geq 2,
\end{equation*}%
and for $\nu >-\frac{1}{2},$
\begin{equation*}
J_{\nu }(z):=\sum\limits_{m=0}^{\infty }(-1)^{m}\left(\frac{z}{2}\right)^{2m+\nu }\left(
m!\ \Gamma (m+\nu +1\right)) ^{-1},\quad z>0,
\end{equation*}%
is the Bessel function of the first kind of order $\nu $.
\end{proposition}
The bounded
non-decreasing function
\begin{equation*}
G(\lambda)=\int_{\{u :\,\left\Vert u \right\Vert <\lambda\}}F(du
),~\lambda\geq 0,
\end{equation*}%
is the spectral distribution function of the homogeneous isotropic random
field. Clearly,
\begin{equation*}
\int_{0}^{\infty }dG(\lambda)=F(\mathbb{R}^{n})=B(0)<\infty\,,
\end{equation*}%
and the inversion formula
\begin{equation}
G(\lambda)=2^{(2-n)/2}\Gamma ^{-1}(n/2)\int_{0}^{\infty }J_{n/2}(\lambda r)\ (\lambda r)^{n/2}%
\frac{B(r)}{r}\ dr  \label{cov5}
\end{equation}%
holds.

For $n=3,$ the formula (\ref{cov4}) becomes
$$
B(r)=\int_{0}^{\infty }\frac{\sin \left( \lambda r\right) \ }{\lambda r}\
G(d\lambda). $$
\begin{definition}\label{def3} We say that $\;\xi (x),~x\in \mathbb{R}^{n}$ is a homogeneous
isotropic random field possessing an absolutely continuous spectrum, if there exists a function $g(\lambda )$ such that
$$
G'(\lambda )=2\pi^{n/2}\Gamma^{-1}\left(n/2\right)\ \lambda^{n-1}g(\lambda
),\quad \lambda^{n-1}g(\lambda )\in L_{1}(\mathbb{R}_{+}).
$$
 The function $g(\lambda )$ is
called the isotropic spectral density function of the field $\xi (x).$
\end{definition}
In this
case, the representation (\ref{cov4}) may be written as
$$
B(r)=2\pi^{n/2}\Gamma^{-1}\left(n/2\right)\ \int_{0}^{\infty }Y_{n}(\lambda r)\
\lambda^{n-1}g(\lambda )\ d\lambda .  $$
\begin{proposition}{\rm(\citealt{leo2})} If $B(r),r\geq 0,$ decays rapidly at infinity, in particular, if
\begin{equation*}
\int_{0}^{\infty }r^{n-1}\left\vert B(r)\right\vert dr<\infty ,
\end{equation*}%
then the function $g(\lambda )$ can be given by the formula
\begin{equation*}
g(\lambda )=\frac{1}{(2\pi )^{\frac{n}{2}}}\int_{0}^{\infty
}J_{(n-2)/2}(\lambda r)\ (\lambda r)^{(2-n)/2}r^{n-1}B(r)\
dr.
\end{equation*}%
\end{proposition}

Let $v(r):=\{x\in R^{n}:\left\Vert x\right\Vert <r\},$ $n\geq 1,$ and $s_{n-1}(r):=\{x\in R^{n}:\left\Vert x\right\Vert =r\},$ $n\geq 2,$ be an open ball and a sphere
of radius $r>0.$ We consider two averaged functionals
$$
\eta (r):=\int_{v(r)}\xi (x)\ dx,\qquad \zeta (r):=\int_{s_{n-1}(r)}\xi (x)\ d\sigma
(x),  $$
where $\sigma (x)$ is the Lebesgue measure on the sphere $s_{n-1}(r).$

Then, the following results hold, see \citealt{yad,leo1,leo2}.

\begin{proposition}\label{pr4}
For $n\ge 1$
\begin{eqnarray*}
b_{n}(r) &:=&{\rm{Var}}\,\eta (r)  =\int_{v(r)}\int_{v(r)}B(\left\Vert x-y\right\Vert )dx\
dy   \\
&=&(2\pi )^{n}r^{2n}\int_{0}^{\infty }J_{n/2}^{2}(\lambda r)\
(\lambda r)^{-n}dG(\lambda ).
\end{eqnarray*}%

If $n\geq 2$ then
\begin{eqnarray*}
l_n(r) &:=&{\rm{Var}}\,\zeta (r)  =\int_{s(r)}\int_{s(r)}B(\left\| x-y\right\| )d\sigma
(x)d\sigma (y)   \\
&=&(2\pi )^nr^{2(n-1)}\int_0^\infty J_{(n-2)/2}^2(\lambda r)\
(\lambda r)^{2-n}dG(\lambda ).
\end{eqnarray*}
\end{proposition}

\section{Regular varying functions}

\begin{definition}\label{def4} A measurable function $L:(0,\infty )\rightarrow (0,\infty )$ is
called slowly varying at infinity if for all $t >0,$%
\begin{equation*}
\lim\limits_{\lambda\to \infty }\frac{L(\lambda t)}{L(\lambda)}=1.
\end{equation*}
\end{definition}

For example, for $k_1>0,\ k_2>0,$ and $k_3 >0,$ the functions
$$
L(\lambda)\equiv 1,\quad L(\lambda)=\left(\log(k_1+\lambda)\right)^{k_2},\quad
L(\lambda)=\log \log (\lambda+k_1),\quad$$
$$L(\lambda)=\frac{\lambda^{k_2k_3}}{\left(k_1+\lambda^{k_2}\right)^{k_3}},\quad L(\lambda)=\exp\left((\log x)^\frac13\cos((\log x)^\frac13)\right),$$
vary slowly at infinity.

Let $\mathcal{L}$ be the class of functions that are slowly varying at
infinity and bounded on each finite interval.

Now we introduce a wider class of $O$-regularly varying functions.
Let
\[L^*(t):=\limsup_{\lambda\to\infty}\frac{L(\lambda t)}{L(\lambda)},\quad
L_*(t):=\liminf_{\lambda\to\infty}\frac{L(\lambda t)}{L(\lambda)},\quad
\lambda>0.\]

\begin{definition}  A positive measurable function $L(\cdot)$ is
$O$-regularly varying (belongs to the class $OR$) if for all $t \ge 1:$
\[0<L_*(t)\le L^*(t)<+\infty\,.\]
\end{definition}

\begin{definition} Let $L(\cdot)$ be a positive measurable function. The infimum of those $\alpha$ for
which there are constants $C=C(\alpha),$ such that
\[ \frac{L(\lambda t)}{L(\lambda)} \le
C\,\{1+o(1)\}t^\alpha\ \ \mbox{uniformly in}\
t\in[1,T]\ \mbox{for all}\ T>1,\ \mbox{as}\ \lambda\to\infty,\]
 is called the upper Matuszewska index and is denoted by $\alpha(L).$

The supremum of those $\beta$ for which there are constants $
D=D(\beta)>0,$ such that
\[ \frac{L(\lambda t)}{L(\lambda)} \ge
D\,\{1+o(1)\}t^\beta\  \ \mbox{uniformly in}\
t\in[1,T]\ \mbox{for all}\ T>1,\ \mbox{as}\ \lambda\to\infty,\]
 is called the lower Matuszewska index and is denoted by $\beta(L).$
\end{definition}

$L(\cdot)\in OR$ if and only if both of its Matuszewska indices $\alpha(L)$ and $\beta(L)$ are
finite.

\begin{definition} By $OR(\beta,\alpha)$ we denote the subclass of $OR$ functions whose Matuszewska indices satisfy $\alpha(L)\le\alpha$ and $\beta\le\beta(L).$
\end{definition}

More details and examples of regular varying functions can be found in \citealt{sen,bingol}.

\section{Classical Abelian and Tauberian  theorems for random fields}

For simplicity we consider only mean-square
continuous homogeneous and isotropic random fields $\xi (x)$, $x\in R^{n}$,
with mean zero and unit variance, i.e. $B(0)=1.$

\citealt{pit} and \citealt{wol} obtained Abelian and Tauberian theorems for  characteristic functions of probability distributions, which corresponds to the one-dimensional case in our notations. The behaviour of the characteristic function in the neighborhood of the origin and the distribution function on infinity was studied.

\citealt{bin}  studied the multidimensional case. He obtained his results in terms of Hankel type transforms (\ref{cov4}).

\begin{theorem}\label{th1}{\rm (\citealt{bin})} Let $0<\gamma <2$, $L(\cdot)\in \mathcal{L},$ then the following two
statements are equivalents:

\begin{enumerate}
  \item[\rm{(a)}] $1-B(r)\sim r^\gamma L\left(\frac{1}{r}\right),$ $r\to 0+;$
  \item[\rm{(b)}] $1-G(\lambda)\sim \frac{2^\gamma
\Gamma \left(\frac{n+\gamma }2\right)}{\Gamma \left(\frac n2\right)\Gamma \left(1-\frac \gamma 2\right)}\cdot \frac{L(\lambda )}{\lambda ^\gamma },$ $\lambda \to \infty .$
\end{enumerate}

If $\gamma =2,$ the statement {\rm{(a)}} is equivalent to
\begin{equation*}
\int\limits_{0}^{\lambda }\mu [ 1-G(\mu )]\ d\mu \sim n\cdot L(\lambda
),\quad \lambda \to +\infty ,
\end{equation*}%
or
\begin{equation*}
\int\limits_{0}^{\lambda }\mu ^{2}dG(\mu )\sim 2n\cdot L(\lambda ),\quad \lambda
\to +\infty.
\end{equation*}%
If $\gamma =0,$  {\rm{(a)}} implies {\rm{(b)}}. Conversely,   the statement  {\rm{(b)}} implies {\rm{(a)}} with $\gamma =0$ if $%
1-G(\lambda )$ is convex for $\lambda $ sufficiently large, but not in
general.
\end{theorem}

\begin{example}
Let $n\ge 3$  and
\begin{equation}\label{G0}
G'(\lambda)= \frac{a^{n-1}\lambda^{n-2}e^{-a\lambda}}{(n-2)!}\,,\quad a> 0\,.
\end{equation}
The corresponding covariance function is (see \citealt[Example 1.2.8]{leo2})
$$B(r)= \frac{a^{n-1}}{(r^2+a^2)^{\frac{n-1}{2}}}\,. $$
It is easy to check that
$$\lim _{r\rightarrow 0} \left( 1-{\frac {{a}^{n-1}}{ \left( {r}^{2}+{a}
^{2} \right)^\frac{n-1}{2}}} \right) {r}^{-2}=
\frac{n-1}{2a^2}.$$
Hence, we obtain $\gamma=2$ and $L\left(\cdot\right)\equiv\frac{n-1}{2a^2}$ in the statement (a) of Theorem~\ref{th1}.

By (\ref{G0}) we get
$$\int\limits_{0}^{\lambda }\mu ^{2}dG(\mu )=\frac{a^{n-1}}{(n-2)!}\int\limits_{0}^{\lambda }\mu^{n}e^{-a\mu}d\mu=
{\frac {{n}^{2}\Gamma  \left( n-1 \right) -n\Gamma  \left( n-1
 \right) -\Gamma  \left( n+1,a\lambda \right) }{{a}^{2}\Gamma  \left( n-1
 \right) }},$$
 where
$\Gamma(c,z) = \int_z^\infty e^{-t}t^{c-1}dt$ is the upper incomplete Gamma function.

Therefore, we  obtain the second statement of Theorem~\ref{th1}:
$$\int\limits_{0}^{\lambda }\mu ^{2}dG(\mu )\sim \frac{n^2-n}{a^2}=2n\cdot L(\lambda ).$$

For example, for $n=3$ and $a=1$ the integral can be calculated explicitly:
$$\int\limits_{0}^{\lambda }\mu ^{2}dG(\mu )=6-6\,{ e^{-\lambda}} \left( 1+\lambda+\frac12\,{\lambda}^{2}+\frac16\,{\lambda}^{3} \right).$$
\end{example}

Various generalizations of Theorem~\ref{th1} were proposed in \citealt{bin1,bin2,bin3,bin4,ino}. The corresponding Tauberian theorems for $\gamma>2$ were obtained by \citealt{soni}. However, results similar to Theorem~\ref{th1} cannot be used to derive Abelian and Tauberian theorems for LRD random fields.

\section{Abelian and Tauberian  theorems for long-range dependent random fields}

Now we are interested in the opposite case of the asymptotic behaviour of $B(r)$ at infinity and $G(\lambda )$ at zero.  Some asymptotic results for Hankel transforms which can be used for this case were obtained by \citealt{soni1}, \citealt{bin0}, \citealt{ino}.

In this section we present theorems of Tauberian and Abelian type for nonintegrable covariance functions of homogeneous isotropic random fields.

Using the inversion formula (\ref{cov5})  we
obtain the representation
$$2^{(n-2)/2}\Gamma(n/2)\cdot\frac{G(\lambda)}{\lambda^{\frac{n-1}{2}}}=\int_{0}^{\infty }\sqrt{\lambda r}\,J_{n/2}(\lambda r)\ r^{\frac{n-3}{2}}B(r)\ dr.$$
By Theorem 1.2 in \citealt{bin0} with $\nu=\frac{n}{2},$ $\alpha=\frac{n}{2}-\rho-\frac{3}{2}$ we get the following result.

\begin{theorem}\label{th2} Let  $L(\cdot)\in \mathcal{L}$ and $r^{\frac{n-3}{2}}B(r)$ be ultimately decreasing to zero at infinity, that is, decreasing for $r$ larger than some number N. Then for $\frac{n-3}{2}<\alpha <n$ the following two statements are equivalent:

\begin{enumerate}
  \item[\rm{(a)}] $r^\alpha B(r)\sim  L\left(r\right),$ $r\to \infty;$
  \item[\rm{(b)}] $G(\lambda )/\lambda^{\alpha}\sim  L\left(\frac{1}{\lambda}\right)/c_{1}(n,\alpha),$ $\lambda
\rightarrow 0+,$
\end{enumerate}
where
\begin{equation}
c_{1}(n,\alpha ):=\frac{2^{\alpha }\Gamma \left( \frac{\alpha}{2}+1\right) \
\Gamma \left( \frac{n}{2}\right) }{\Gamma \left( \frac{n-\alpha }{2}\right) }\,.%
\   \label{const4}
\end{equation}
\end{theorem}

\begin{remark} In the one-dimensional case the constant~(\ref{const4}) can be simplified to the form $
c_{1}(1,\alpha )=\Gamma (1+\alpha )\cos \left(\frac{\alpha \pi }{2}\right).
$
\end{remark}
 The case when $\alpha > n$ can also  be obtained from corresponding results by \citealt{ino} but it has a more complicated version of the statement (b) with additional terms and needs moore delicate conditions.

For some values of $\alpha$ it is even possible to obtain similar results without assumptions about ultimately decreasing behaviour of covariance functions. In \citealt{ole0} and \citealt{leo2} it was shown that:

\begin{theorem}\label{th3} For $0<\alpha <n$ the statement {\rm{(a)}} in Theorem~{\rm\ref{th2}} implies the statement {\rm{(b)}.} For  $0<\alpha <(n-3)/2,$ $n\geq 4,$ the statement {\rm{(b)}} implies the statement {\rm{(a)}.}
\end{theorem}

Now we present an example of asymptotic behaviours in Theorem~\ref{th3}.

\begin{example}
Suppose that $n=9,$ $\alpha=2,$ and
    $$G(\lambda)= \left\{
     \begin{array}{ll}
       \lambda^2, & 0\le  \lambda \le a; \\
       a^2, & \lambda>a.
     \end{array}
     \right.$$
Then the corresponding covariance function is
$$B(r)= {\frac {14 \left( {a}^{5}{r}^{5}+45\,\cos \left( ar
 \right) ra-45\,\sin \left( ar \right) +15\,\sin \left( ar \right) {r}
^{2}{a}^{2} \right) }{{a}^{5}{r}^{7}}}
. $$
It is obvious that the statements of Theorem~\ref{th3} are satisfied.

For $a=2$ plots of  $B(r)$  and its normalized transformation  are shown in Fig.~1.

\begin{center}
\begin{minipage}{6cm}
{\psfig{figure=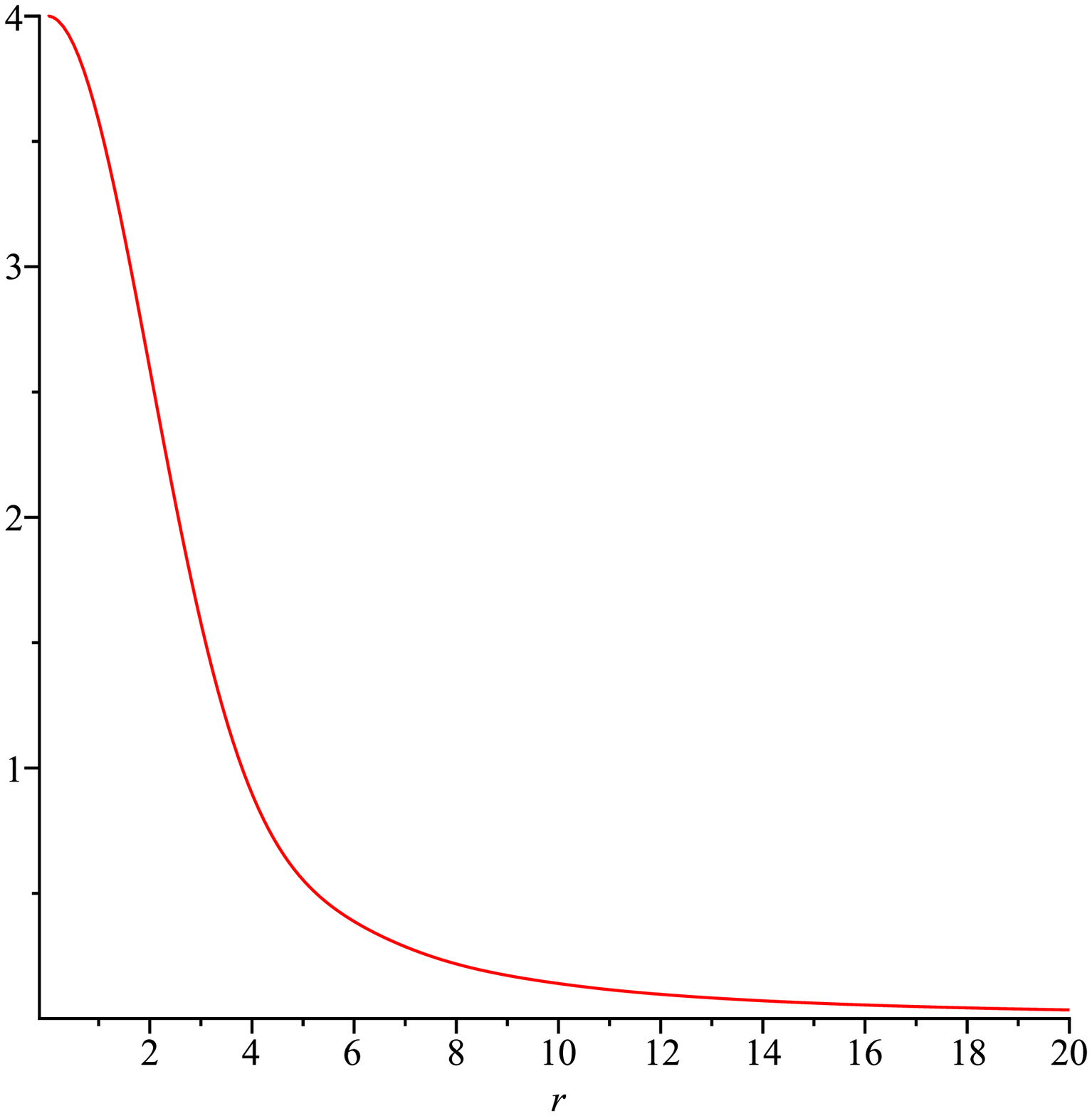 ,width=5cm}} \hspace{0.5cm}{  {\textbf{Fig.\,1a}}\ \   \small Plot of  $B(r)$} \end{minipage}
\begin{minipage}{6cm}
{\psfig{figure=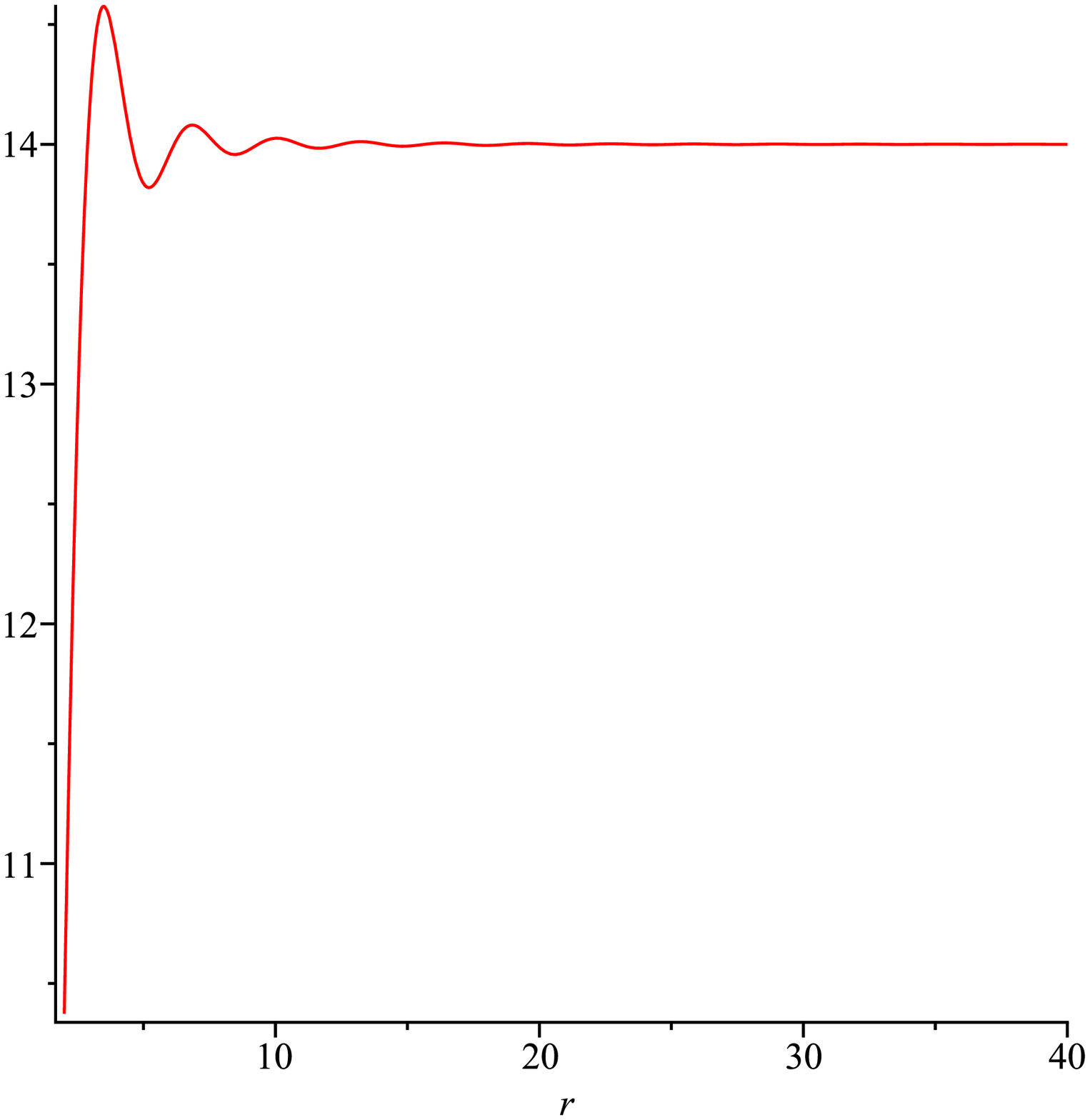,width=5cm}} \hspace{0.5cm}{ {\textbf{Fig.\,1b}} \ \ \small Plot of  $r^2 B(r)$} \end{minipage}
\end{center}
\end{example}

Under some additional assumptions it is possible to prove similar results in terms of isotropic spectral densities too.
For example,  the next result follows from Theorem~\ref{th3}, see more general Theorem 7 in \citealt{ole1}.
\begin{theorem}\label{th4} Suppose that there exists an isotropic spectral density $g(\lambda),$ $\lambda \in [0,\infty),$  such that $g(\lambda)$ is decreasing for all $\lambda\in (0, \varepsilon],$  $\varepsilon > 0.$  Then for $0<\alpha <n$ the statement {\rm{(a)}} in Theorem~{\rm\ref{th3}} implies
\begin{enumerate}
\item[\rm{(b${}^\prime$)}]
$\lambda^{n-\alpha}\,g(\lambda)\sim  L\left(\frac{1}{\lambda}\right)/c_{2}(n,\alpha ),$ $\lambda \rightarrow 0, $
where $c_{2}(n,\alpha ) :=\frac{2^{\alpha }\pi
^{n/2}\Gamma \left( \frac{\alpha }{2}\right) }{\Gamma \left( \frac{n-\alpha }{2}\right) }.$
\end{enumerate}
\end{theorem}

\begin{remark} For $n=1$ we can simplify the constant to $c_2(1,\alpha )=2\Gamma (\alpha ) \cos \left(\frac{\alpha \pi }2\right).
$
\end{remark}

\begin{remark} Note that the class of characteristic functions of symmetric distributions in $\mathbb{R}^n,$
coincides with the class of correlation functions of homogeneous and isotropic random fields. Thus, we can use this duality fact to construct covariance functions, in which properties of densities of multivariate symmetric distributions become properties of spectral densities of random fields.
\end{remark}

\begin{example}\label{ex2}
(\emph{Bessel covariance function})
Let us consider the covariance function
\begin{equation}
B\left( r\right) =\frac{1}{\left( 1+r^{2}\right) ^{\kappa/2}},\ \kappa>0,\ r\ge 0.  \label{bessel1}
\end{equation}%
The function (\ref{bessel1}) is the characteristic function of the
multivariate Bessel distribution, see \citealt[p.69]{fan}. In geostatistics the function (\ref{bessel1}) is known as the Cauchy covariance function, see \citealt{gne}.

Due to {\rm(\ref{bessel1})} the statement (a) of Theorem~\ref{th2} holds.

 The corresponding isotropic spectral density, see {\rm\citealt[p. 292]{don}}, is
\begin{equation}
g\left( \lambda \right) =\left( \pi ^{\frac{n}{2}}2^{\frac{%
n+\kappa  -2}{2}}\Gamma \left( \frac{\kappa }{2}\right) \right) ^{-1}K_{\frac{%
n-\kappa}{2}}\left(  \lambda  \right)
\lambda ^{\frac{\kappa-n}{2}},~\lambda\ge 0,
\label{bessel2}
\end{equation}%
where $K_{\nu}(z)$ is the modified Bessel function of the second kind.

The asymptotic behaviour of  $K_{{}\nu }\left(
z\right) $ at zero is
\begin{equation}
K_{{}\nu }(z)\sim \frac12 \Gamma(\nu)\left(\frac{z}{2} \right)^{-\nu},\ \nu>0.\label{besselk}
\end{equation}%
Substituting (\ref{besselk}) into (\ref{bessel2}) we see that the statement~\rm{(b${}^\prime$)} of Theorem~\ref{th4} holds too.

For $n=6$ and $\kappa=4$ plots of  $g\left( \lambda\right)$  and its normalized transformation  are shown in Fig.~2. In Fig.~2b the spectral density is plotted in log-log coordinates. It is clear that the logarithmic transform of the spectral density approaches a strait line at negative infinity. The negative slope $\alpha-n=-2$ implies $\alpha=4.$

\

\noindent\begin{minipage}{6cm}
 {\psfig{figure=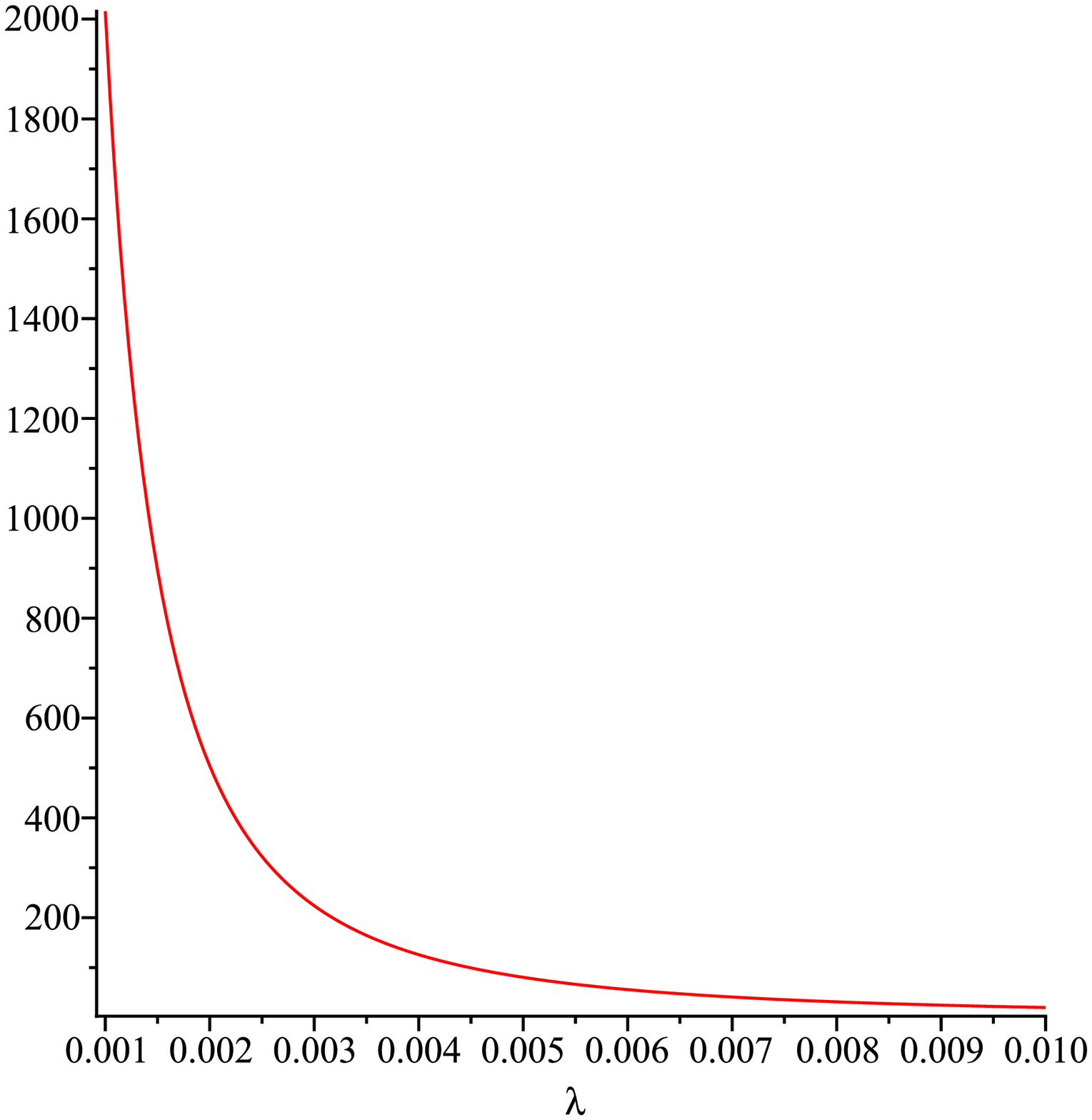,width=5cm}} \hspace{1cm}{\textbf{Fig.\,2a}\ \
 \small Plot of $g\left( \lambda\right)$} \end{minipage}\quad\ \
\begin{minipage}{6cm}
{\psfig{figure=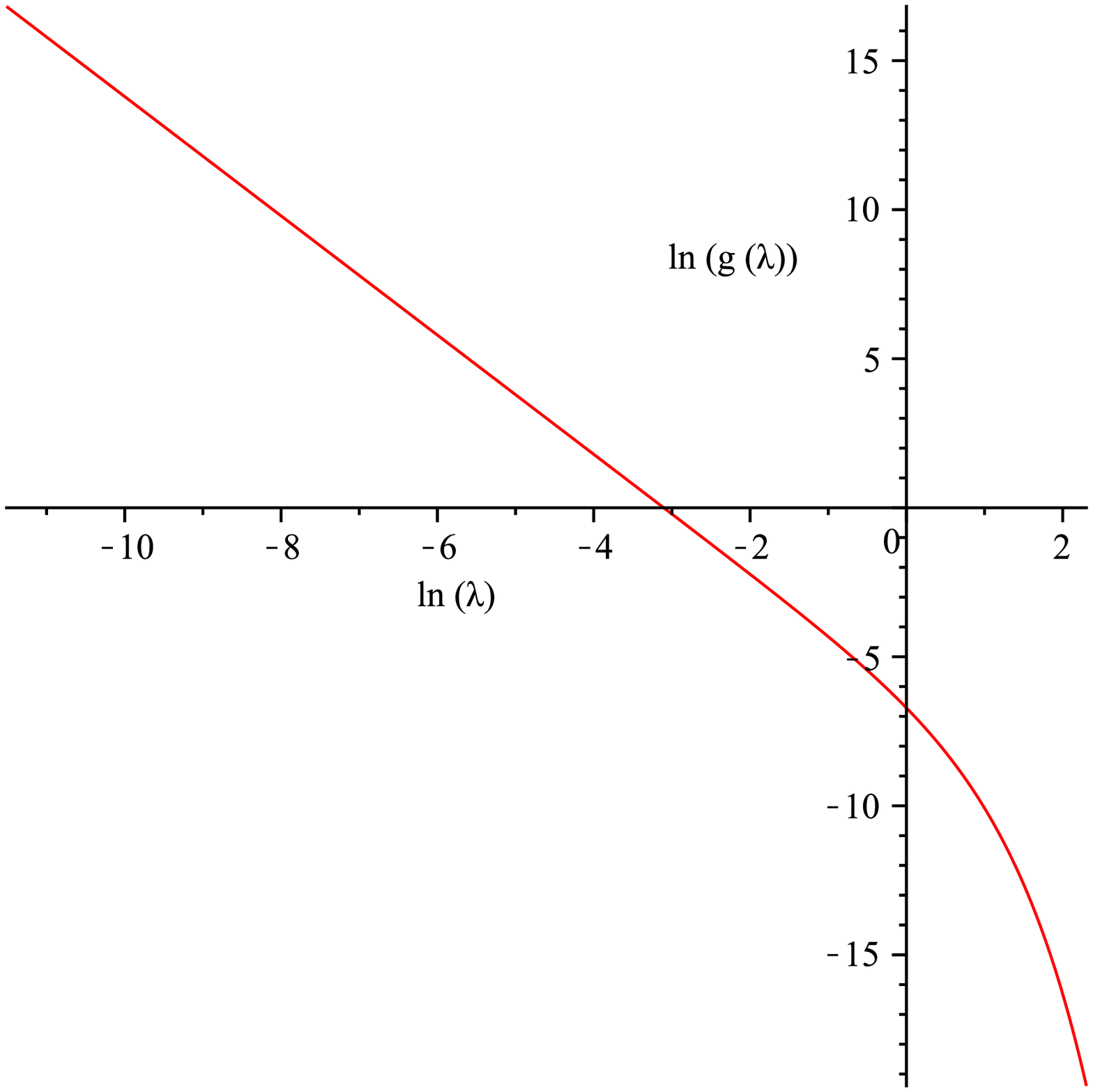,width=5cm}} \hspace{1cm}{\textbf{Fig.\,2b}\ \  \small Log-log plot of $g\left( \lambda\right)$} \end{minipage}
\begin{center}\begin{minipage}{6cm}
 {\psfig{figure=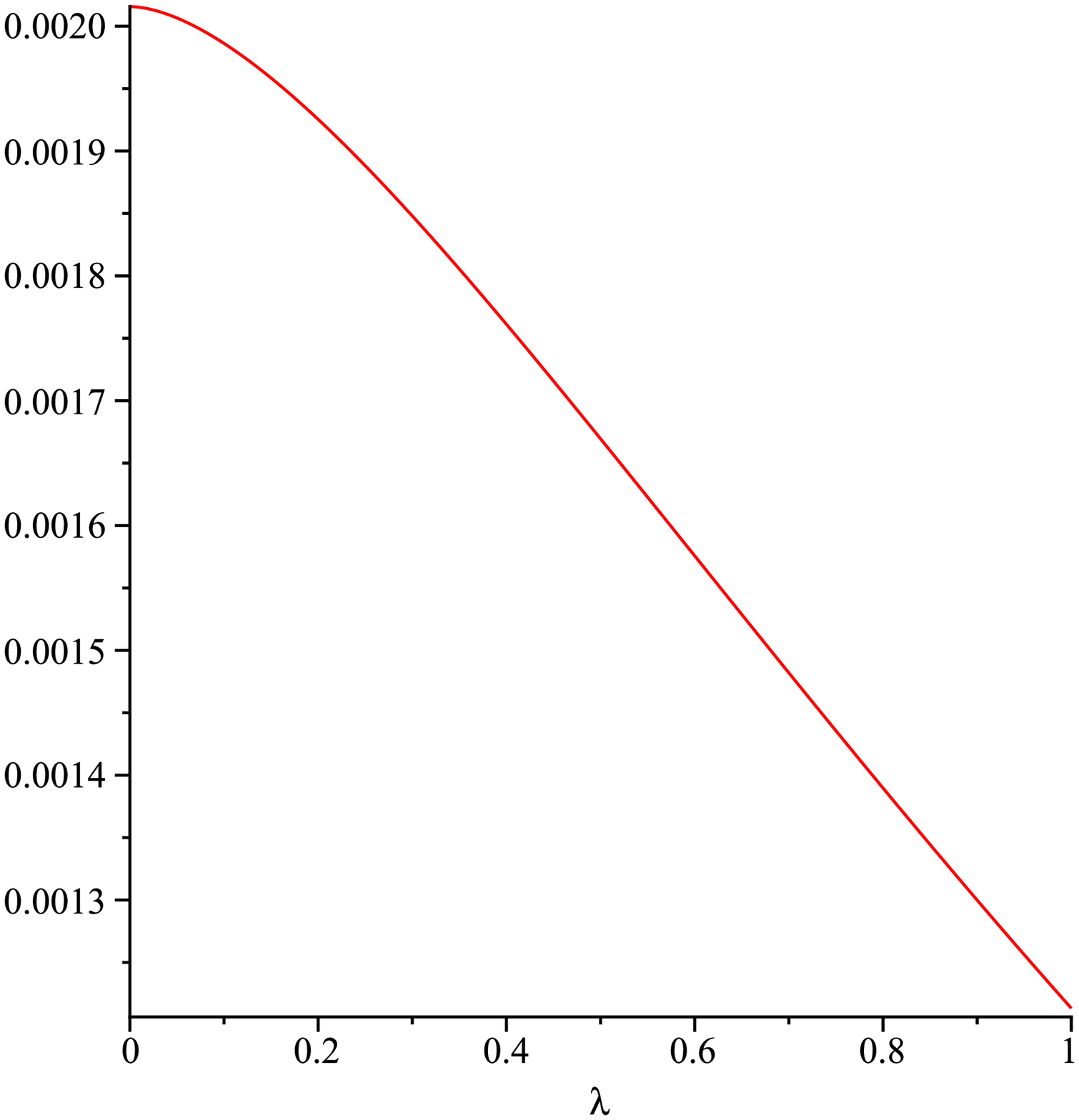,width=5cm}} \hspace{1cm}{\textbf{Fig.\,2c}\ \  \small Plot of
$\lambda^{2}\,g(\lambda)$} \end{minipage}
 \end{center}

\end{example}

\begin{remark}
The covariance function~{\rm(\ref{bessel1}) } and its isotropic spectral density~{\rm(\ref{bessel2})} have the remarkable property that all the convolutions of the spectral
density can be computed:
\begin{equation*}
g_{{}\kappa}^{\ast 1}\left( u \right):=g(\lambda),\quad g_{{}\kappa }^{\ast m}\left( u \right):=\int_{\mathbb{R}%
^{n}}g_{{}\kappa}^{\ast \left( m-1\right) }\left( u^{\prime }\right)
g_{{}\kappa}^{\ast 1}\left( u -u^{\prime }\right) du^{\prime
},\ \ m\ge 2,
\end{equation*}%
where $\lambda:=\left\Vert u\right\Vert,$ $u\in \mathbb{R}^{n}.$

Then, for Bessel covariance functions, every $m\geq 1,$ and $r:=\left\Vert x\right\Vert$ we get
$$
\left( B\left( r\right) \right) ^{m} =\frac{1}{\left(
1+\left\Vert x\right\Vert ^{2}\right) ^{\frac{\kappa m}{2}}} =\int_{\mathbb{R}^{m}}e^{i\left\langle u,x\right\rangle
}g_{{}\kappa}^{\ast m}\left( u \right) du,
$$
from which we obtain the elegant formula
$$
g_{{}\kappa}^{\ast m}\left(u\right) =g_{{}m\kappa}^{\ast 1}\left(u
\right) ,\quad m\geq 2.$$
This formula can be used to construct other Bessel-type examples of Abelian and Tauberian theorems for random fields.
\end{remark}

\begin{example}
(\emph{Linnik covariance function}) Let us consider the covariance function
\begin{equation}
B_\kappa\left( r\right) =\frac{1}{1+r^{{}\kappa}%
},\ \ 0<\kappa \le 2,\ \ r\ge 0.  \label{lin1}
\end{equation}%
The function {\rm(\ref{lin1})} is known as the characteristic function of the
Linnik distribution, see {\rm\citealt{lin,fan,and,ost,erd}}. For $\kappa\in (0,2)$ the isotropic spectral
density of the Linnik covariance function can be represented as
$$
g_\kappa\left( \lambda \right) =\frac{\sin \left( \frac{\pi \kappa}{2}\right) }{2^{%
\frac{n}{2}-1}\pi ^{\frac{n}{2}+1}\lambda ^{\frac{n}{2%
}-1}} \int_{0}^{\infty }K_{^{\frac{n}{2}-1}}\left( \lambda
 u\right) \frac{u^{\frac{n}{2}+\kappa}du}{\left\vert
1+u^{{}\kappa}e^{i\pi \kappa/2}\right\vert ^{2}},~\lambda \ge 0.$$

In \citealt{erd, lim} the generalized Linnik
distribution was considered. Its covariance function takes
the form
\begin{equation*}
B_{{}\kappa,\nu }\left( r\right) =\left( 1+r^{{}\kappa}\right)
^{-\nu },\quad \kappa \in \left( 0,2\right] ,\ \nu >0,\ r\ge 0.
\end{equation*}%
For $\kappa \in \left( 0,2\right),$ $\nu >0,$ the corresponding isotropic spectral density has the form
\begin{equation*}
g_{{}\kappa,\nu }\left( \lambda \right)
=\frac{\lambda^{1-\frac{n}{2}}}{2^{%
\frac{n}{2}-1}\pi ^{\frac{n}{2}+1}} \int_{0}^{\infty }K_{^{\frac{n}{2}-1}}\left( \lambda
 u\right) \frac{\sin(\nu\, \mbox{arg}(1+u^{{}\kappa}e^{i\pi \kappa/2}))}{\left\vert
1+u^{{}\kappa}e^{i\pi \kappa/2}\right\vert ^{\nu}}\,u^{\frac{n}{2}}du,~\lambda \ge 0.
\end{equation*}
By Corollary~{\rm 3.10} in {\rm\citealt{lim} } for $\kappa \nu <n:$
\begin{equation*}
g_{{}\kappa,\nu }\left( \lambda \right)\sim \frac{\Gamma\left(\frac{n-\kappa\nu}{2}\right)}{2^{\kappa\nu}\pi^{\frac{n}{2}}\Gamma\left(\frac{\kappa\nu}{2}\right)}\lambda^{\kappa\nu-n},\quad  \lambda \to 0,
\end{equation*}
and we get the statements {\rm (a)} and {\rm (b${}^\prime$)}.

For example, for $n=3,$ $\kappa=1$ and $\nu=2$ we obtain
$$g_{1,2}\left(\lambda\right)=\frac {\sin \left(\lambda\right)\left(2\,{\it Ci} \left(\lambda\right)+2\lambda\,{\it Si} \left(\lambda\right) -\lambda\pi\right)+ \cos \left(\lambda \right)\left(\pi-2\,{
\it Si} \left(\lambda\right)  +2\lambda\,{\it Ci} \left(\lambda \right)\right)}{4\lambda{\pi }^{2}},$$
where
 $$ Ci(x)=\gamma+\ln(x)+\int_0^x\frac{\cos(t)-1}{t}dt,  \;\;\;
   Si(x)=\int_0^x\frac {\sin(t)}{t}dt,$$
and $\gamma$  is the Euler's  constant.

Plots of  $g_{1,2}\left( \lambda\right)$  and its normalized transformation  are shown in Fig.~3. In Fig.~3b the spectral density is plotted in log-log coordinates. From the slope value $\alpha-n=-1$ at negative infinity we get $\alpha=2,$ which is in agreement with the statement {\rm (b${}^\prime$)}.

\

\noindent\begin{minipage}{6cm}
 {\psfig{figure=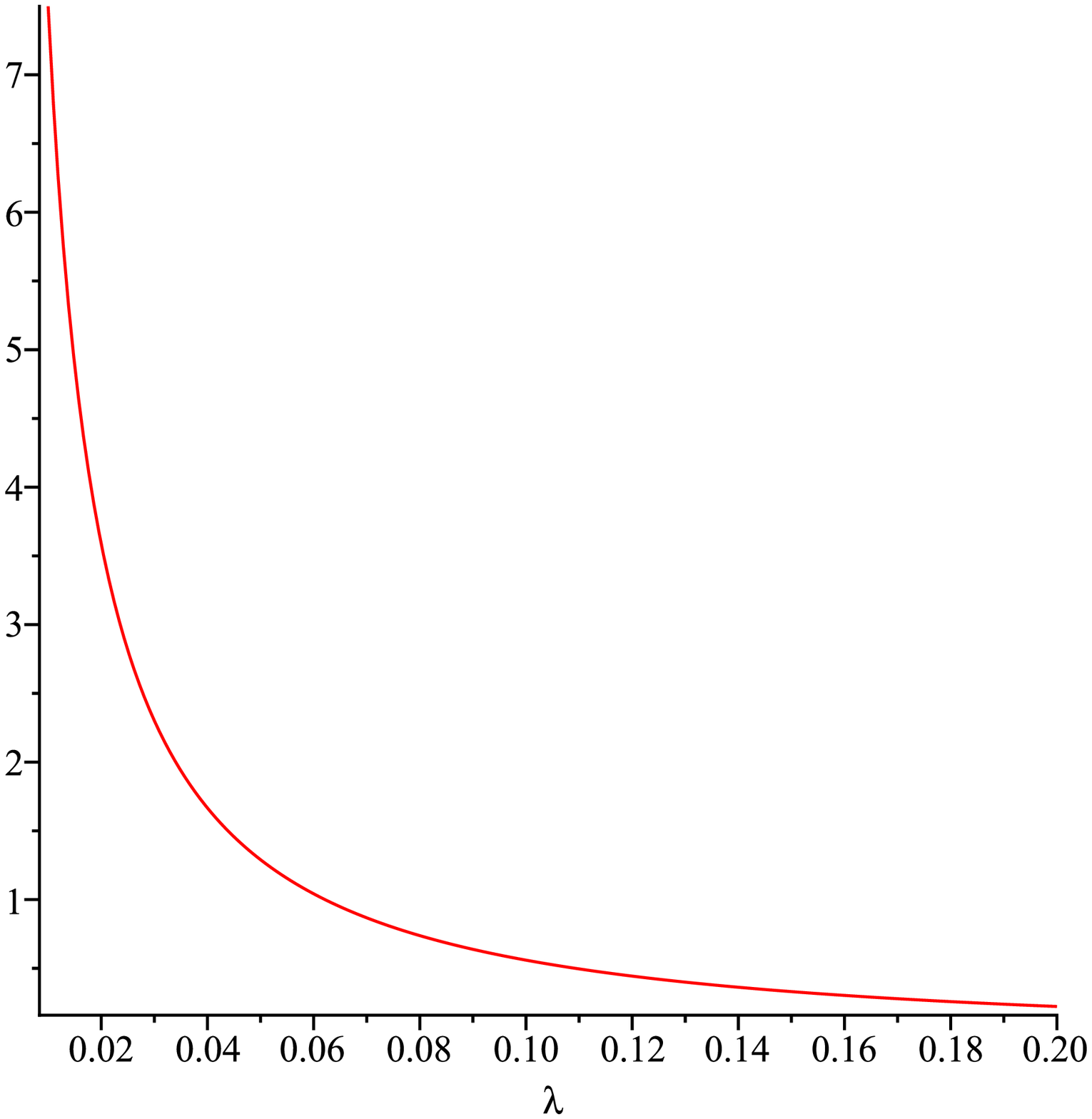,width=5cm}} \hspace{1cm}{\textbf{Fig.\,3a}\ \
 \small Plot of $g_{1,2}\left( \lambda\right)$} \end{minipage}\quad\ \
\begin{minipage}{6cm}
{\psfig{figure=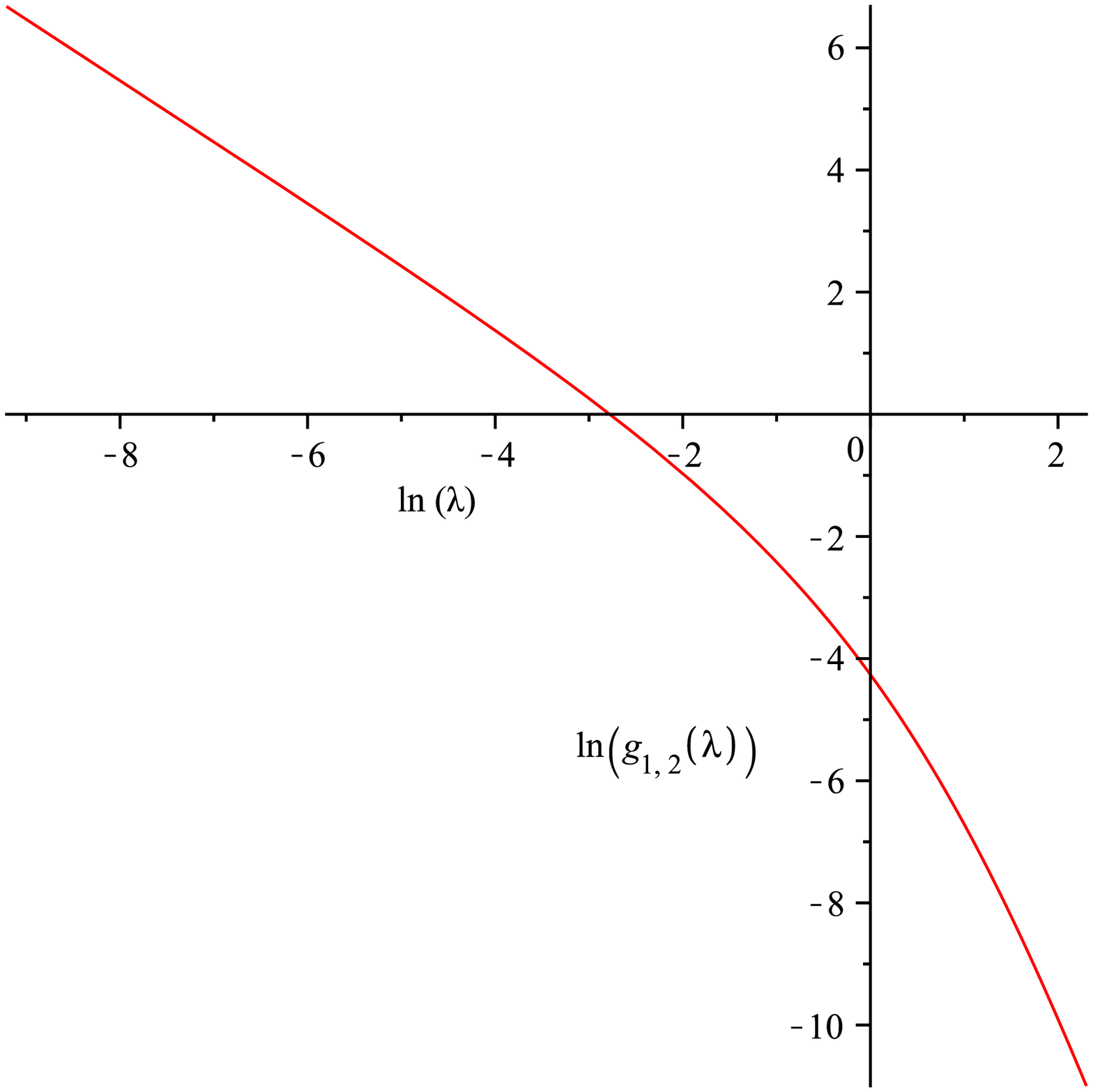,width=5cm}} \hspace{1cm}{\textbf{Fig.\,3b}\ \  \small Log-log plot of $g_{1,2}\left( \lambda\right)$} \end{minipage}
\begin{center}\begin{minipage}{6cm}
{\psfig{figure=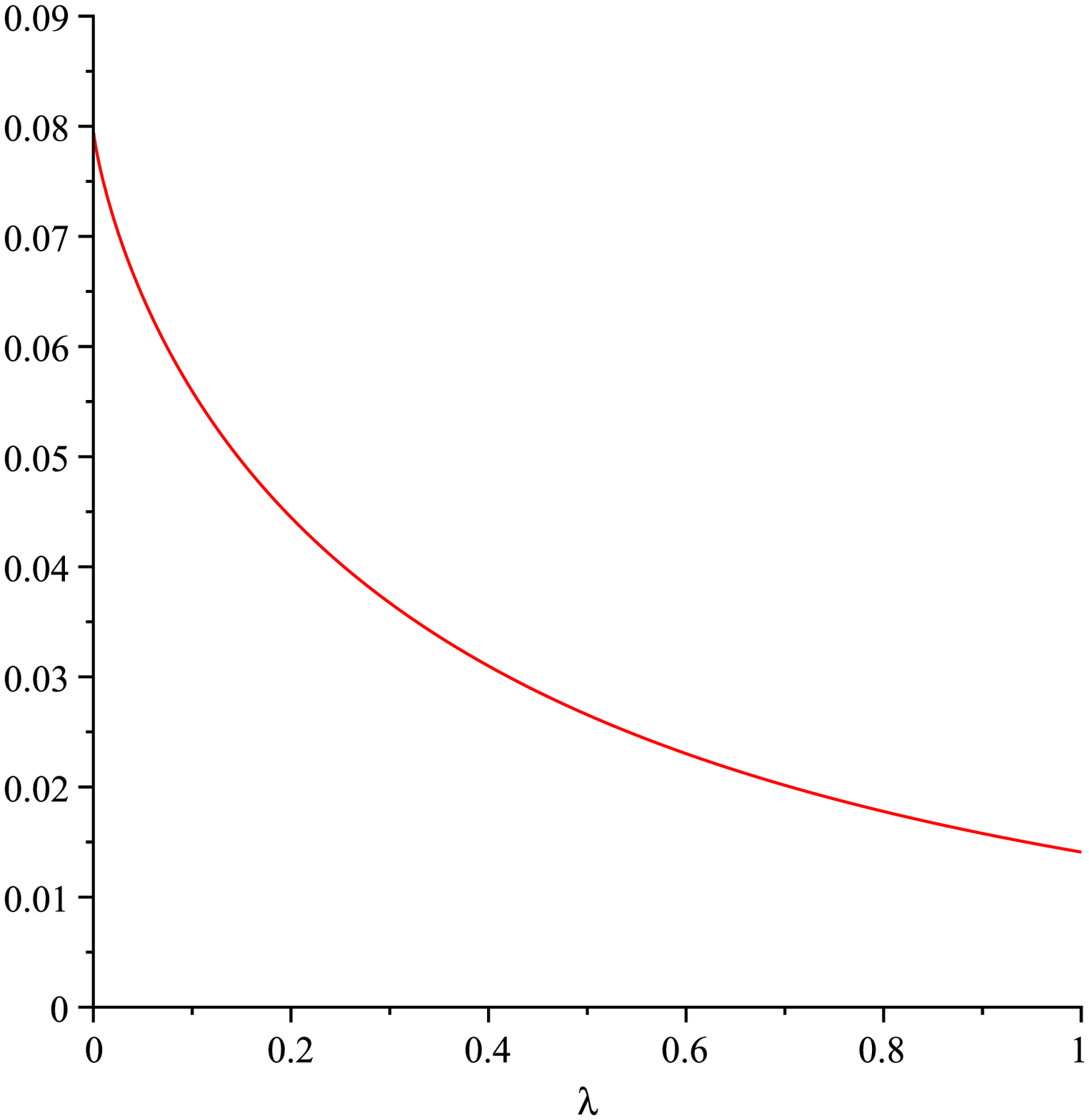,width=5cm}} \hspace{1cm}{\textbf{Fig.\,3c}\ \  \small Plot of
$\lambda\,g_{1,2}(\lambda)$} \end{minipage}
 \end{center}
\end{example}

\section{Abelian and Tauberian  theorems for variances of averaged functionals}

In subsequent sections we explain how Abelian and Tauberian results can be extended to accommodate a more complex asymptotic behaviour. The averaged functionals $b_{n}(\cdot)$ and $l_{n}(\cdot)$ introduced in Proposition~\ref{pr4} are very useful in these studies.

The following example shows that the statements (a) and (b) in Theorem~\ref{th3} are not equivalent for $(n-3)/2<\alpha.$ The equivalence is only true under some additional conditions, see, for example, Theorem~\ref{th2}.

\begin{example}\label{ex4} Let $n=3,$ $\alpha=2,$
    $$G(\lambda)= \left\{
     \begin{array}{ll}
         \lambda^{2}, & 0\le  \lambda \le a; \\
      a^2, & \lambda>a.
     \end{array}
     \right.$$
Then the corresponding covariance function is
$$B(r)= \frac {2 \left(1-\cos \left( ar \right)  \right) }{{r}^{2}},$$
and the variance of averaged functional has the form
$$b_3(r) =4\, \pi^{2}a^{-4} \left( 2\,{r}^{4}{a}^{4}-2\,{r}^{2}{a}^{2}+2\,
\sin \left( 2\,ra \right) ra-1+\cos \left( 2\,ra \right)  \right)
.$$

For $a=1$ plots of  $B(r),$  $b_3(r),$  and their normalized transformations  are shown in Fig.~4.

\begin{center}
\begin{minipage}{5cm}
 {\psfig{figure=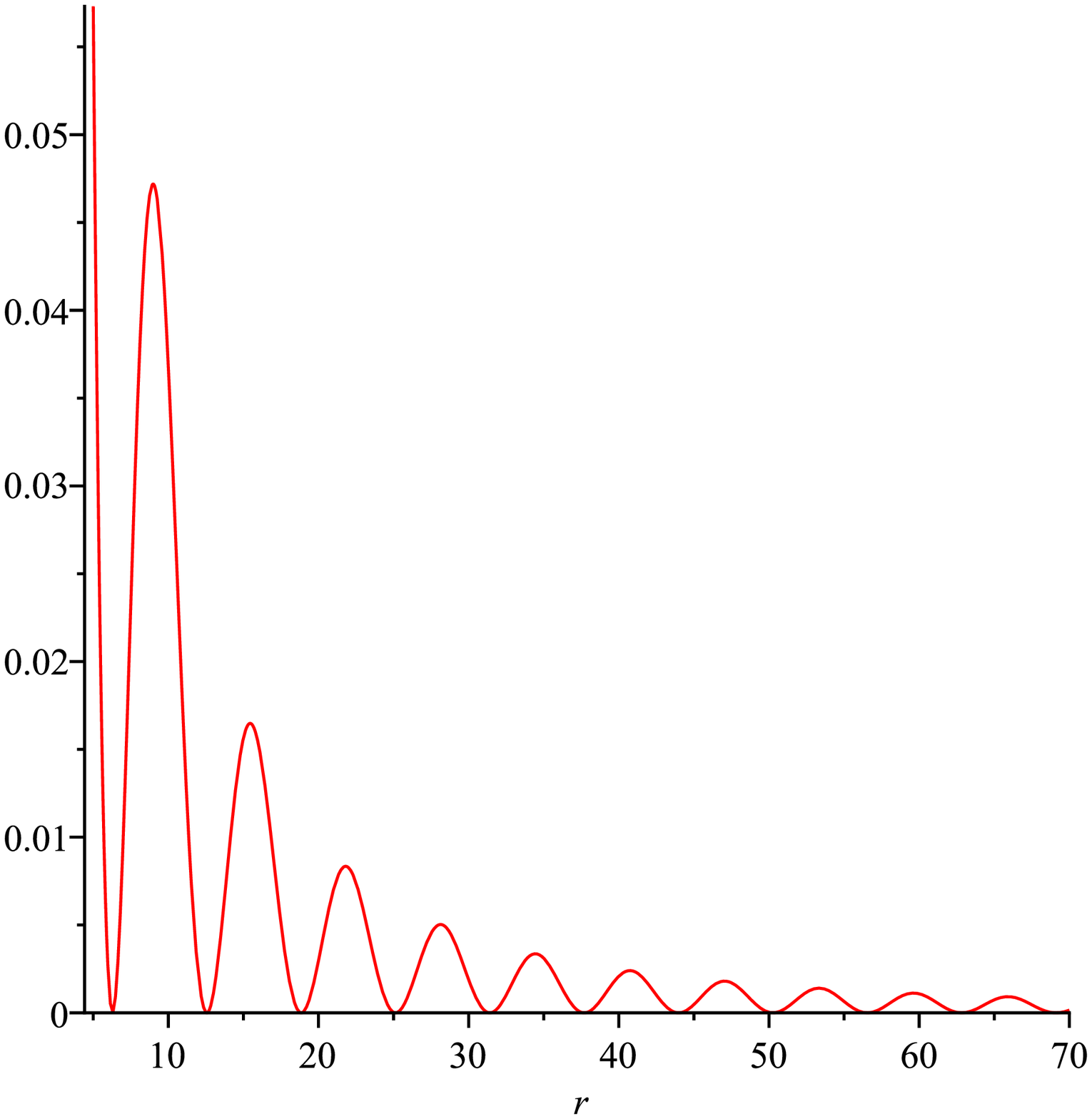 ,width=5cm}} \hspace{0.5cm}{\textbf{Fig.\,4a}\ \ \small Plot of  $B(r)$} \end{minipage}\quad\ \
\begin{minipage}{5cm}
{\psfig{figure=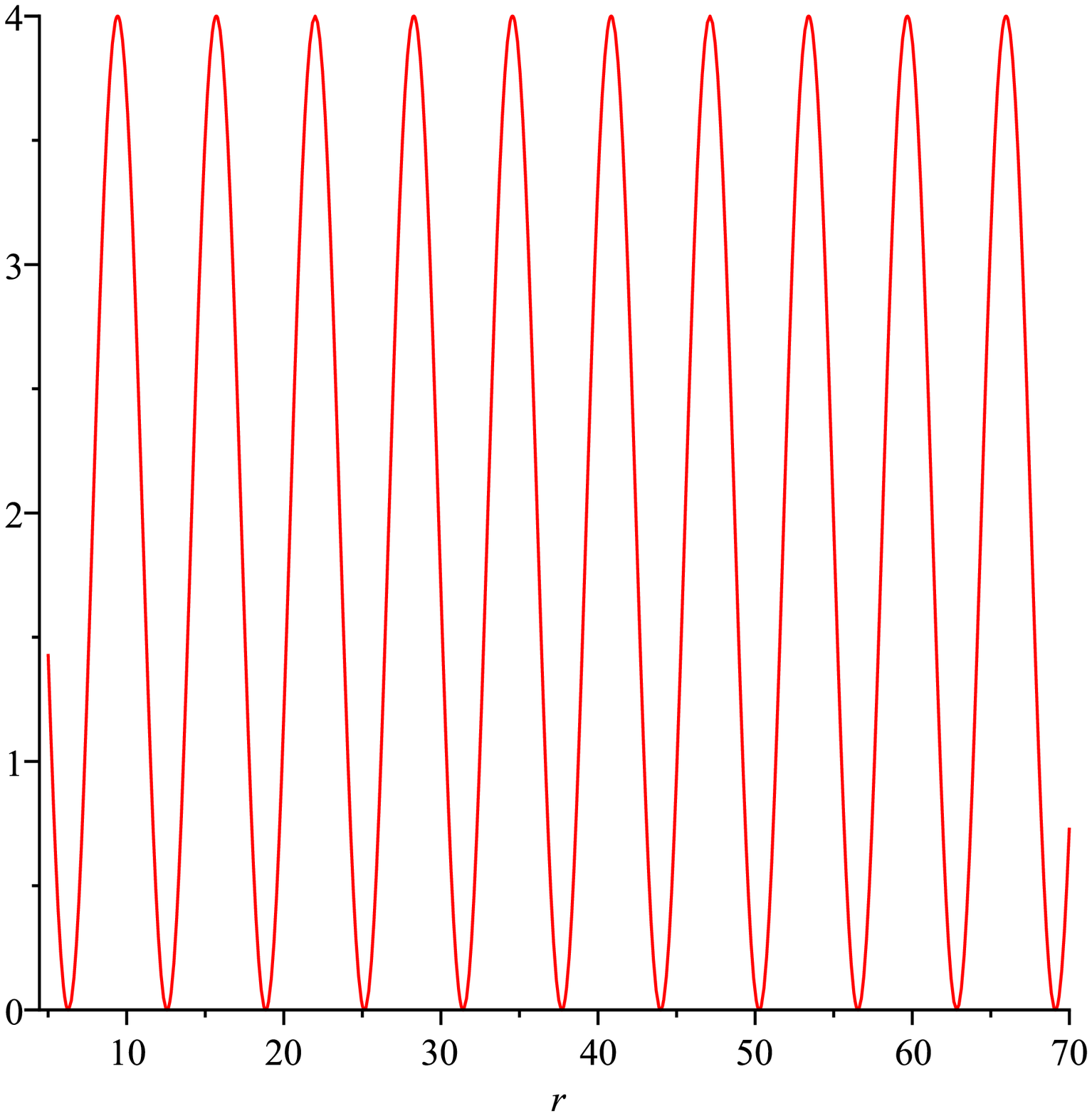,width=5cm}} \hspace{0.5cm}{\textbf{Fig.\,4b}\ \ \small Plot of  $r^2 B(r)$} \end{minipage}
\end{center}
\begin{center}
\begin{minipage}{5cm}
 {\psfig{figure=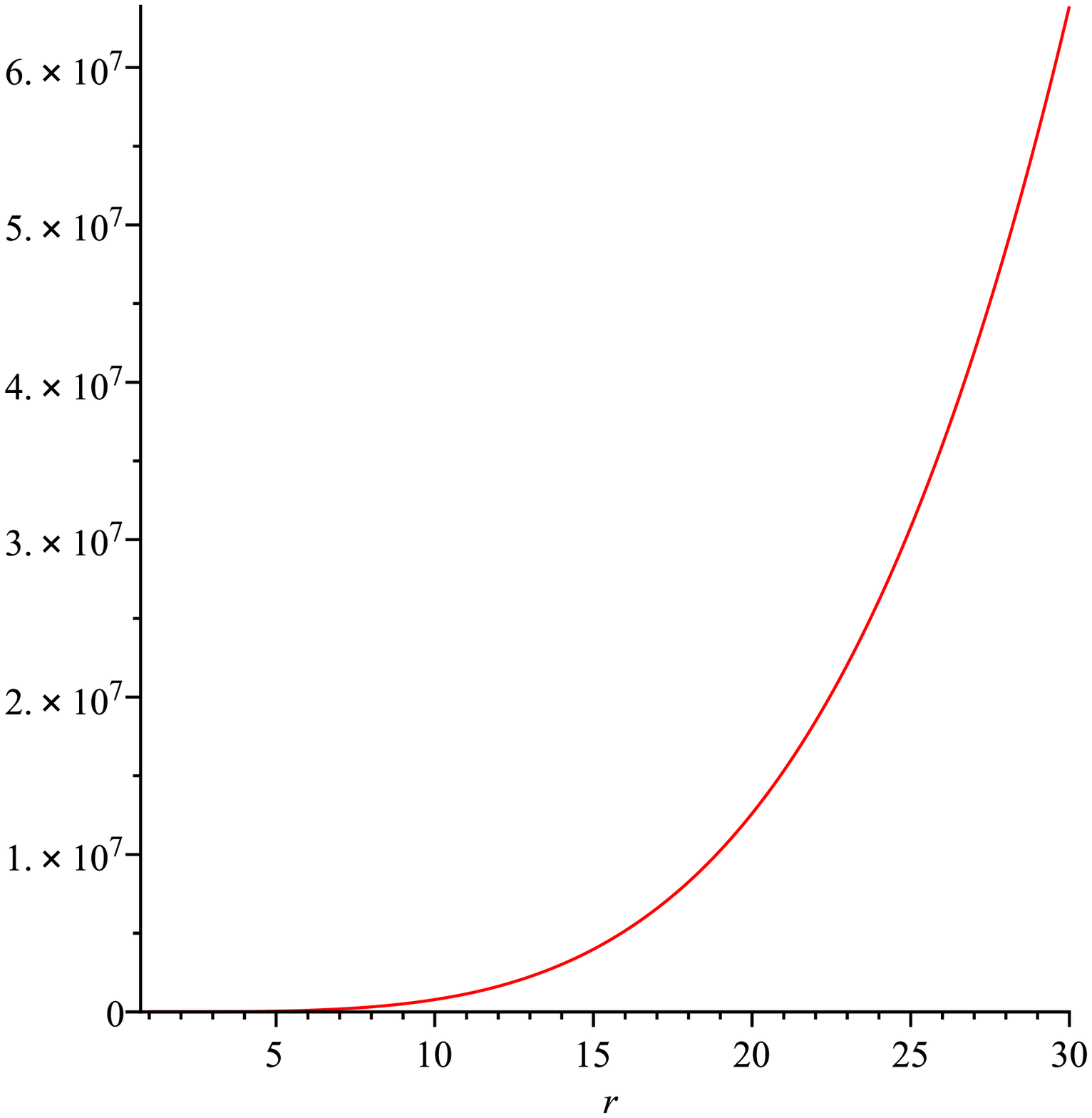 ,width=5cm}} \hspace{0.5cm}{\textbf{Fig.\,4c}\ \  \small  Plot of  $b_3(r)$} \end{minipage}\quad\ \
\begin{minipage}{5cm}
{\psfig{figure=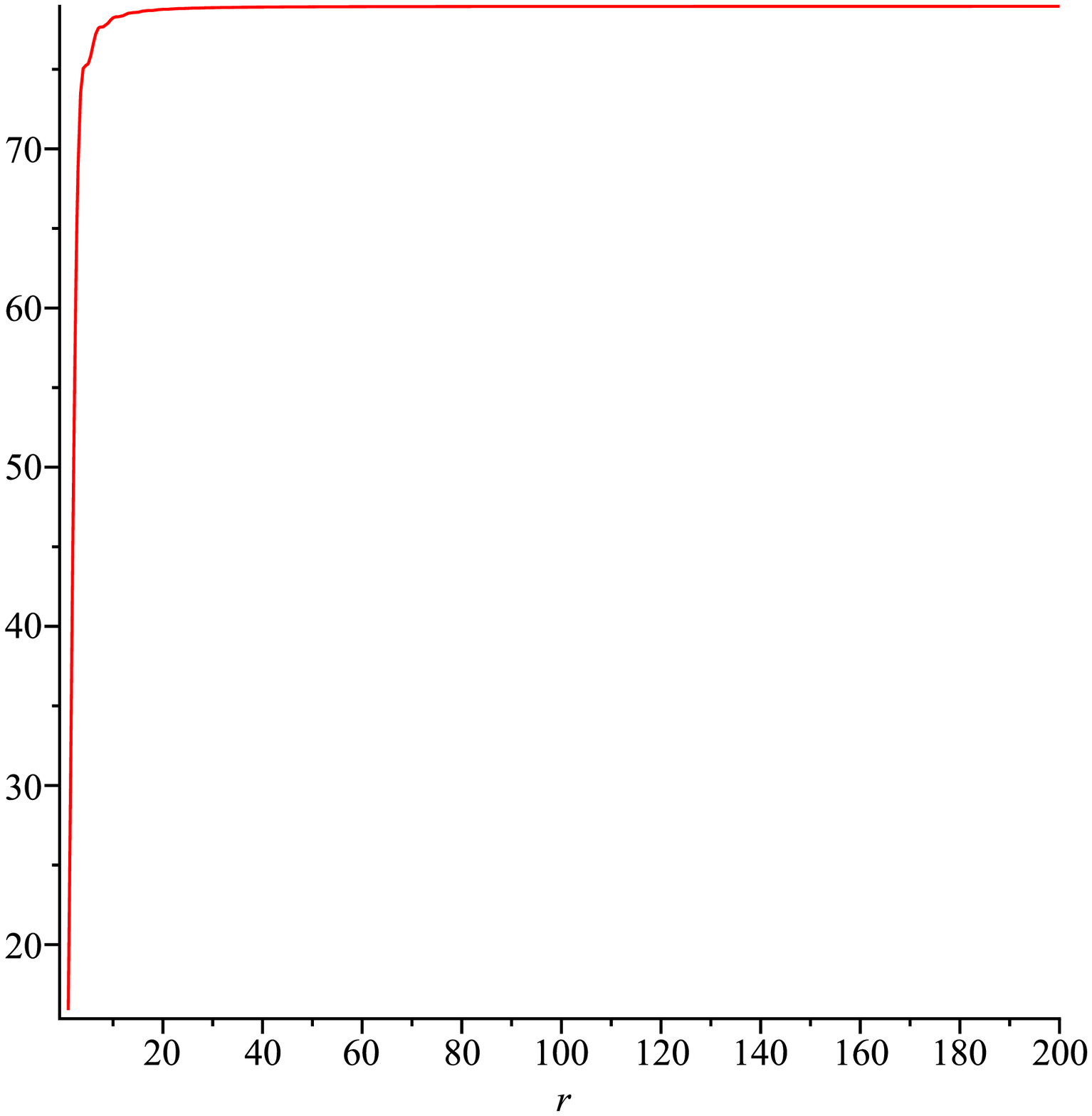,width=5cm}} \hspace{0.5cm}{\textbf{Fig.\,4d}\ \ \small Plot of  $r^{-4} b_3(r)$} \end{minipage}
\end{center}
\end{example}

\begin{example}\label{ex41}
Suppose that $n=3,$ $\alpha=2,$ and
   $$g(\lambda)= \left\{
     \begin{array}{ll}
        \frac{2+\cos(\lambda)}{\lambda}, & 0\le \lambda \le 1; \\
        \lambda^{-1}, & 1< \lambda \le 2; \\
       0, & \;\; \textrm{otherwise}.
     \end{array}
     \right.$$
Then the correlation function of the field is
 $$B(r) = \frac {4\pi}{{r}^{2} \left( {r}^{2}-1 \right) } \left(  \left( 4-\cos \left( r \right) -\cos
 \left( r \right) \cos \left( 1 \right) -2\, \left( \cos \left( r
 \right)  \right) ^{2} \right) {r}^{2}\right.$$
 $$\left.-r\sin \left( r \right) \sin
 \left( 1 \right) -3+\cos \left( r \right) +2\, \left( \cos \left( r
 \right)  \right) ^{2} \right).$$
The formula for the variance of averaged functional $b_3(r)$ is omitted, but it has closed-form representation as well.

Plots of  $B(r),$  $b_3(r),$  and their normalized transformations  are displayed in Fig.~5. Some empirical covariance functions in spatial statistics demonstrate similar oscillatory behaviour.

\begin{center}
\begin{minipage}{5cm}
 {\psfig{figure=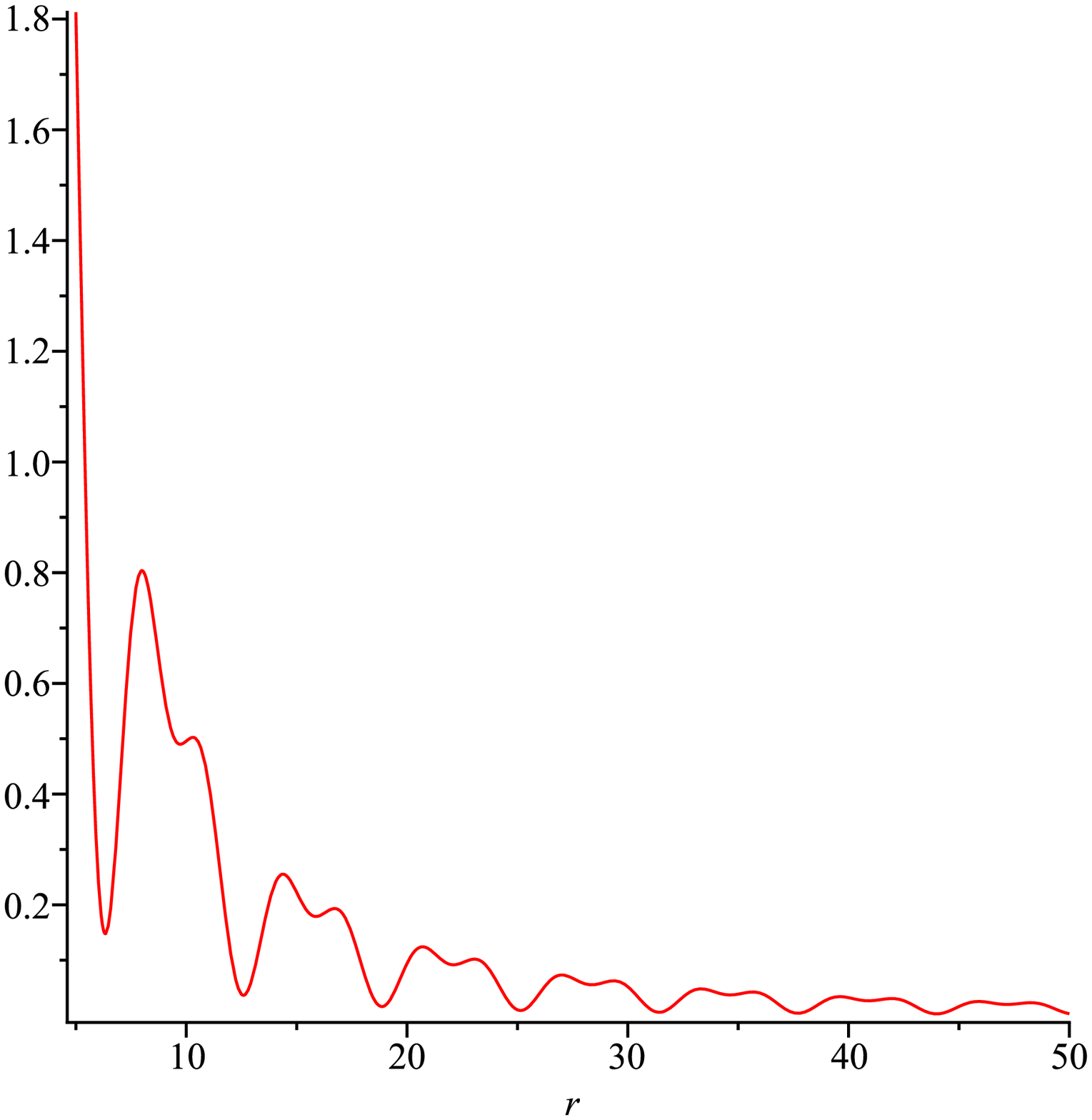 ,width=5cm}} \hspace{0.5cm}{\textbf{Fig.\,5a}\ \ \small Plot of  $B(r)$} \end{minipage}\quad\ \
\begin{minipage}{5cm}
{\psfig{figure=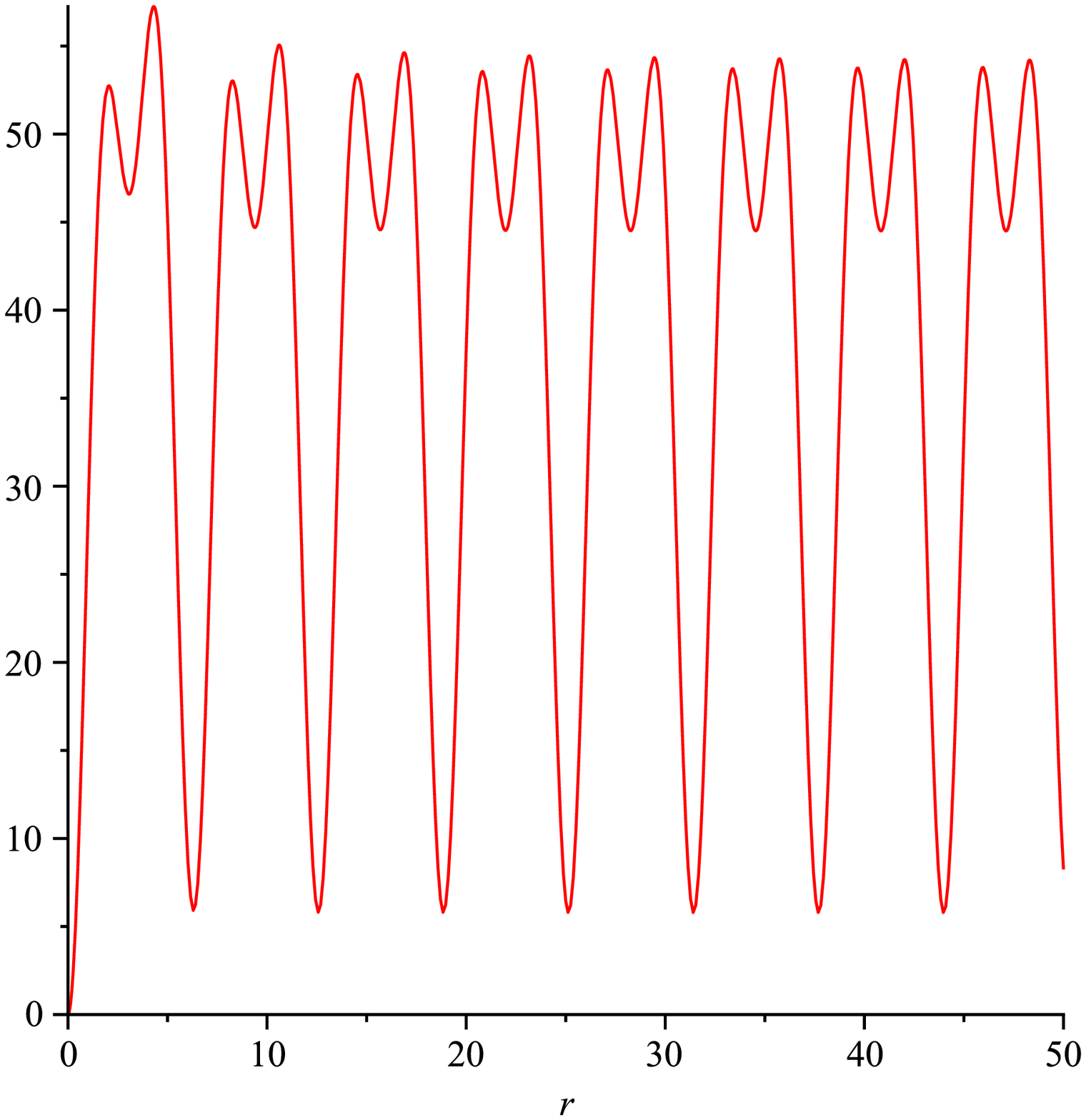,width=5cm}} \hspace{0.5cm}{\textbf{Fig.\,5b}\ \ \small Plot of  $r^2 B(r)$} \end{minipage}
\end{center}
\begin{center}
\begin{minipage}{5cm}
 {\psfig{figure=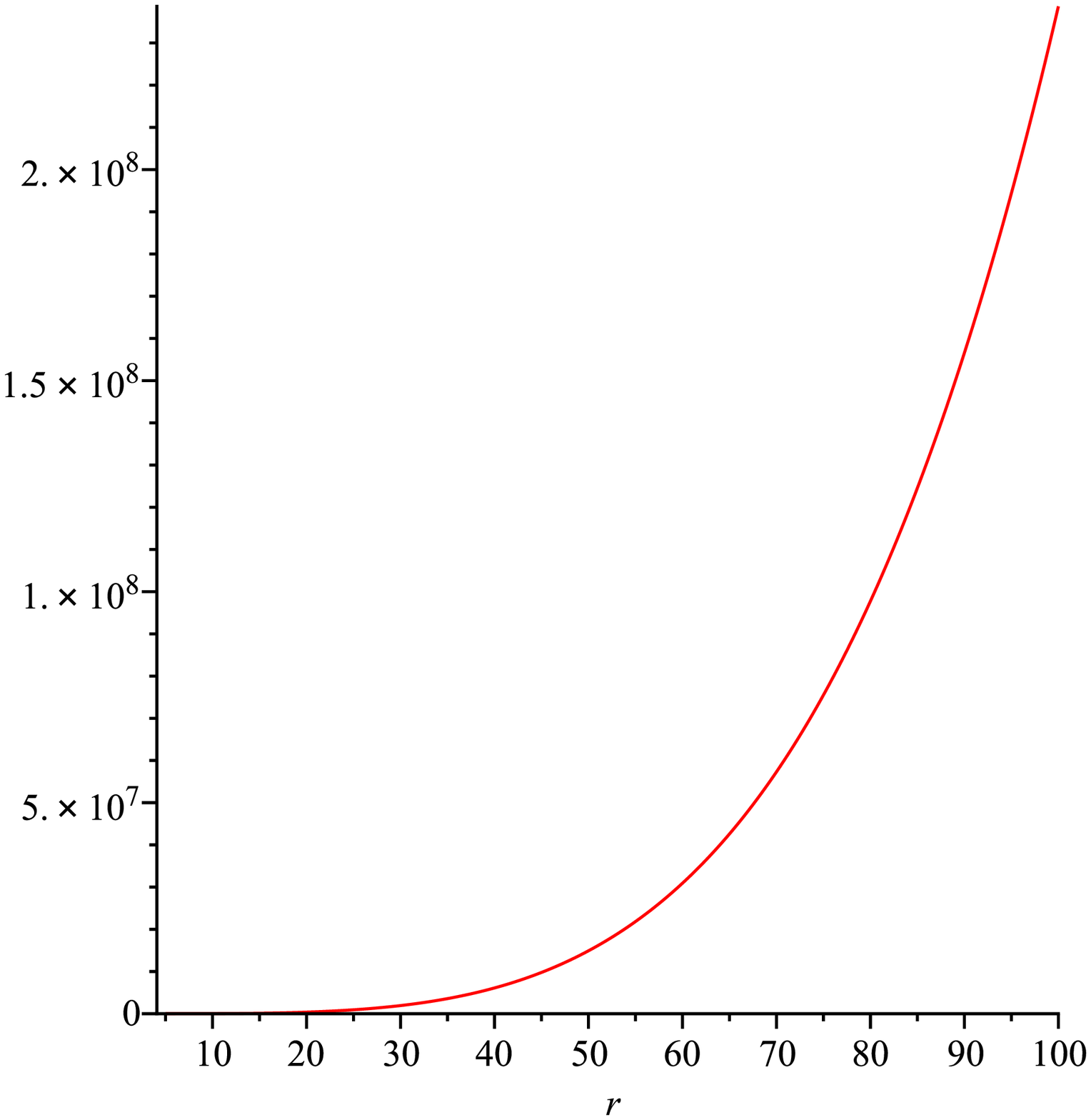 ,width=5cm}} \hspace{0.5cm}{\textbf{Fig.\,5c}\ \  \small  Plot of  $b_3(r)$} \end{minipage}\quad\ \
\begin{minipage}{5cm}
{\psfig{figure=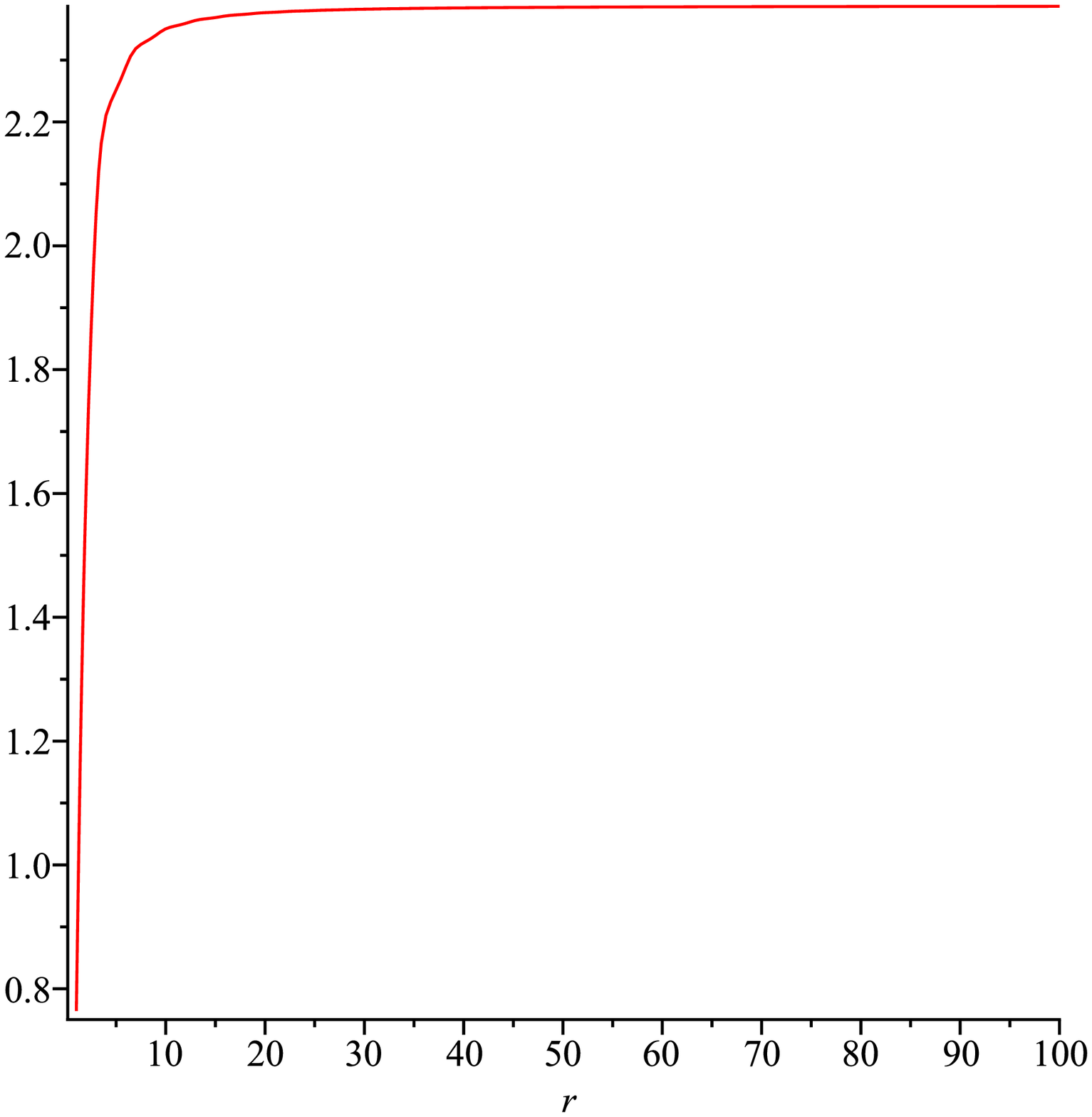,width=5cm}} \hspace{0.5cm}{\textbf{Fig.\,5d}\ \ \small Plot of  $r^{-4} b_3(r)$} \end{minipage}
\end{center}
The functions $G(\lambda)$ and $g(\lambda)$ in Examples~\ref{ex4} and \ref{ex41} satisfy the statement (b) of Theorem~\ref{th2}  and the statement (b${}^\prime$) of Theorem~\ref{th4}  respectively, but the corresponding  function $r^2 B(r)$ is not regularly varying when $r\to +\infty.$ However, the function $b_3(r)$ exhibits regular varying behaviour at infinity.
\end{example}

Examples~\ref{ex4} and \ref{ex41} show that it might be reasonable to use the variances of averaged functionals instead of $B(r)$ to obtain Abelian and Tauberian theorems.  Indeed, one-dimensional results of such type were proved to be true by \citealt{lau}. Some special cases of such multidimensional Abelian theorems were given in \citealt{yad,leo1,leol1,leol2}. The general results for random fields were presented by \citealt{ole0, ole96, leo2}:

\begin{theorem} Let $0<\alpha <n-1,$ $n\geq 2,$ $L(\cdot)\in
\mathcal{L}$. The following two statements are equivalent:
\begin{enumerate}
  \item[\rm{(a)}] $l_n(r)/r^{2n-\alpha -2}\sim L(r),\quad r\rightarrow \infty\,;$
\item[\rm{(b)}]
$G(\lambda )/\lambda ^\alpha \sim L(\frac 1\lambda)/c_3(n,\alpha ),\quad \lambda \rightarrow 0+\ ,$
\end{enumerate}
where
$$
c_3(n,\alpha ):=\frac{\alpha \pi ^{n}2^{\alpha +1}\Gamma (n-\alpha -1)\
\Gamma \left( \frac{\alpha }{2}\right) }{\Gamma ^{2}\left( \frac{n-\alpha }{2%
}\right) \ \Gamma \left( n-1-\frac{\alpha }{2}\right) }.\   $$
\end{theorem}

\begin{theorem} Let $0<\alpha <n,$ $L(\cdot)\in \mathcal{L}$. The
following two statement are equivalent:
\begin{enumerate}
  \item[\rm{(a)}] $b_{n}(r)/r^{2n-\alpha }\sim L(r),\quad
r\rightarrow \infty\,;$
\item[\rm{(b)}] $G(\lambda )/\lambda ^\alpha \sim L(\frac 1\lambda )/c_4(n,\alpha),\quad \lambda \rightarrow 0+\ ,$
\end{enumerate}
where
$$
c_4(n,\alpha ):=\frac{\alpha \ \pi ^{n}2^{\alpha }\Gamma (n-\alpha -1)\
\Gamma \left( \frac{\alpha }{2}\right) }{\Gamma ^{2}\left( \frac{n-\alpha +2%
}{2}\right) \ \Gamma \left( \frac{2n-\alpha +2}{2}\right) } . $$
\end{theorem}

\section{$O$-regularly varying asymptotic behaviour}

Now we shall  give finer results about $O$-regularly varying asymptotic behaviours.

The following example shows that the statements (a) and (b${}^\prime$) in Theorem~\ref{th4} are not equivalent ((b${}^\prime$) does not follow from (a) in  the general case).

\begin{example}\label{ex5}

Suppose that  $n=3$ and $$ g(\lambda)= \frac{e^{-50 \lambda} + \sin^2\left(\frac 2 \lambda \right) }{(2\pi)^{\frac 32}\,\lambda^{\frac 52}},\quad \lambda>0. $$

It is obvious that $\lambda^2 g(\lambda)\in L_1(\mathbb{R}_{+})$  and $g(1/\cdot)\in OR.$ However, there does not exist $\alpha$ for which $\lambda^{3-\alpha}\,g(\lambda)\sim  L\left(\frac{1}{\lambda}\right)/c_{2}(n,\alpha ),$ $\lambda \rightarrow 0.$ To see that one can also look at the plot of the spectral density in Fig.~6a and the log-log plot in Fig.~6b. The logarithmic transform of the spectral density does not approach a strait line at negative infinity.

The corresponding covariance function is
$$B(r) = \frac{ 16\,\sqrt {\sqrt {2500+{r}^{2}}-50}+8\,\sqrt {r}+{{\rm e}^{-4\,\sqrt {r}}}-\sin \left( 4\,
\sqrt {r} \right) -\cos \left( 4\,\sqrt {r}
 \right) }{8r}
.$$

The integral in the representation of $l_3(r)$ can be also explicitly evaluated and we get
$$
 l_3(r)  = \frac{{\pi }^{2}{r}^{2}}{24} \left(128006-(3+12\,\sqrt {2r})\sin
 \left( 4\,\sqrt {2r} \right) +2^8\,\sqrt {2}\,{r}^{3/2}-12\,\sqrt {2r}\,{{\rm e}^{-4\,\sqrt {2r}}}\right.$$
 $$+(12\,\sqrt {2r}-3)\cos \left( 4\,\sqrt {2r}
 \right)-3\,{{\rm e}^{-4\,
\sqrt {2r}}}  \left. +2^9\,\sqrt {50+2\,\sqrt {{r}^{2}+5^4}}\left(\sqrt {{r}^{2}+5^4
} -50\right)\right)
. $$
We see that $r^{\frac 12} B(r)$ and $r^{-\frac 72} l_3(r)$ belong to $\cal{L}.$

Plots of  $B(r),$  $l_3(r),$  and their normalized transformations  are shown in Fig.~6.

\begin{center}
\begin{minipage}{5cm}
 {\psfig{figure=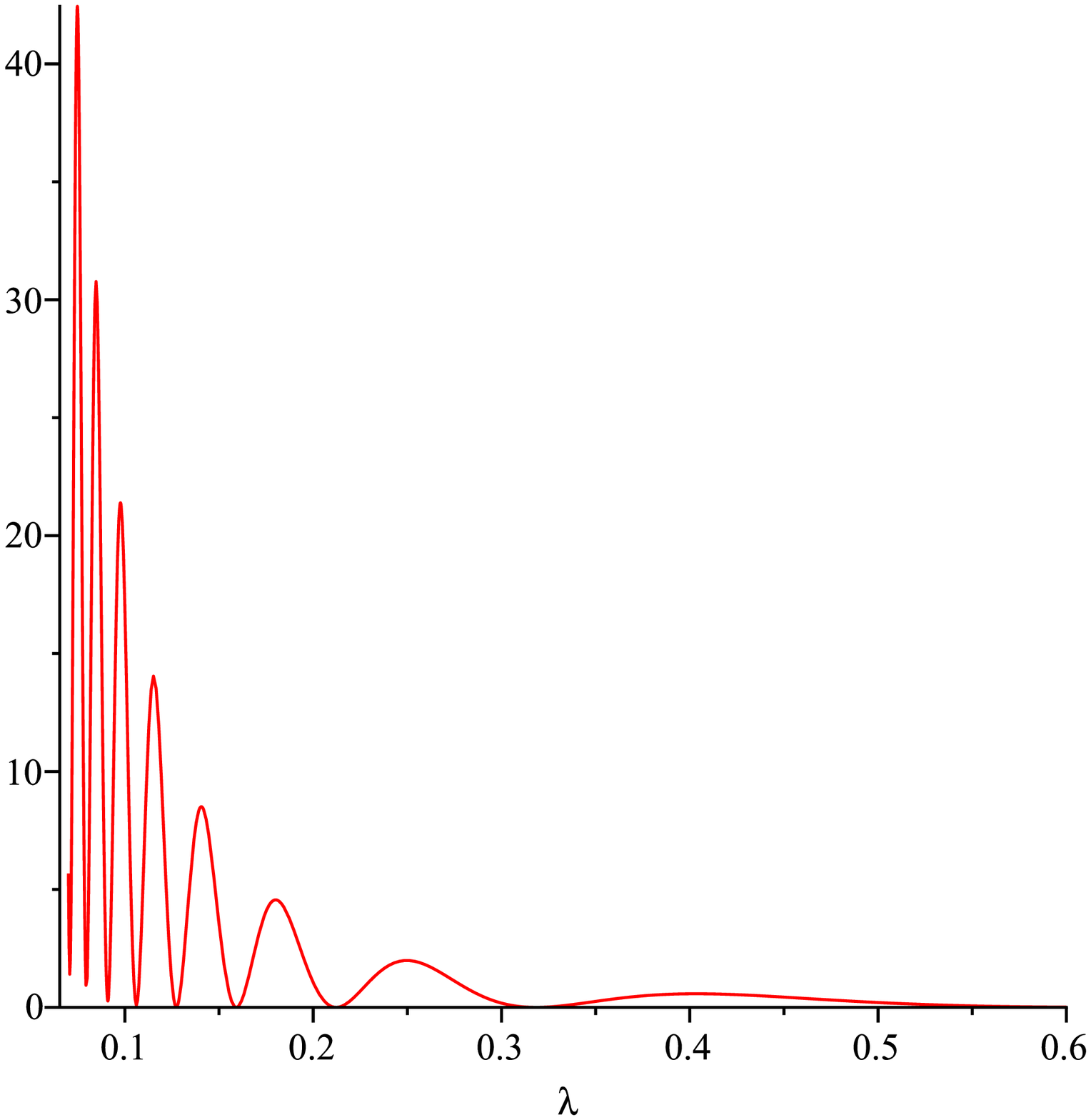 ,width=5cm}} \hspace{0.5cm}{\textbf{Fig.\,6a}\ \ \small Plot of $g\left( \lambda\right)$} \end{minipage}\quad\ \
\begin{minipage}{5cm}
{\psfig{figure=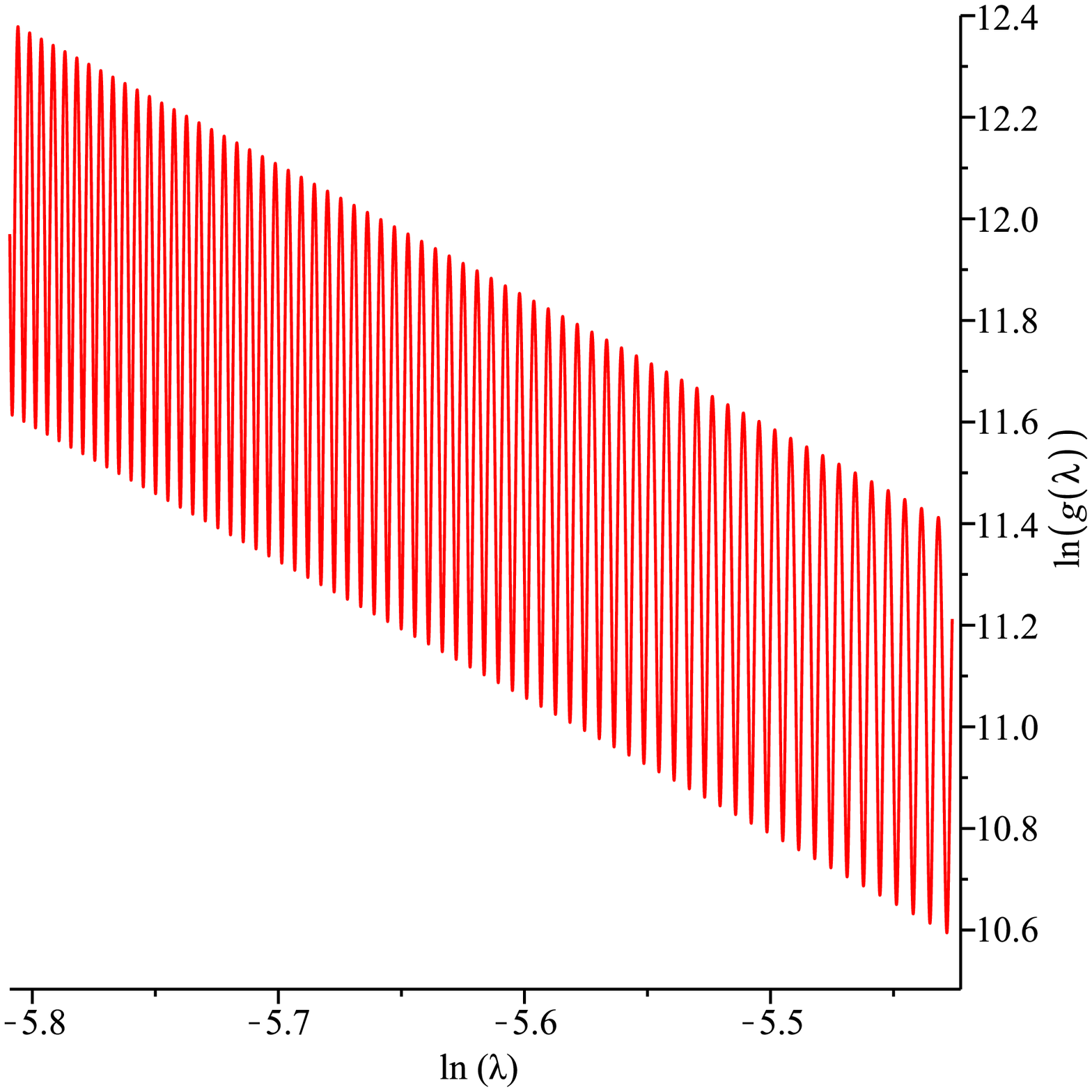,width=5cm}} \hspace{0.5cm}{\textbf{Fig.\,6b}\ \ \small Log-log plot of $g\left( \lambda\right)$} \end{minipage}
\end{center}

\begin{center}
\begin{minipage}{5cm}
 {\psfig{figure=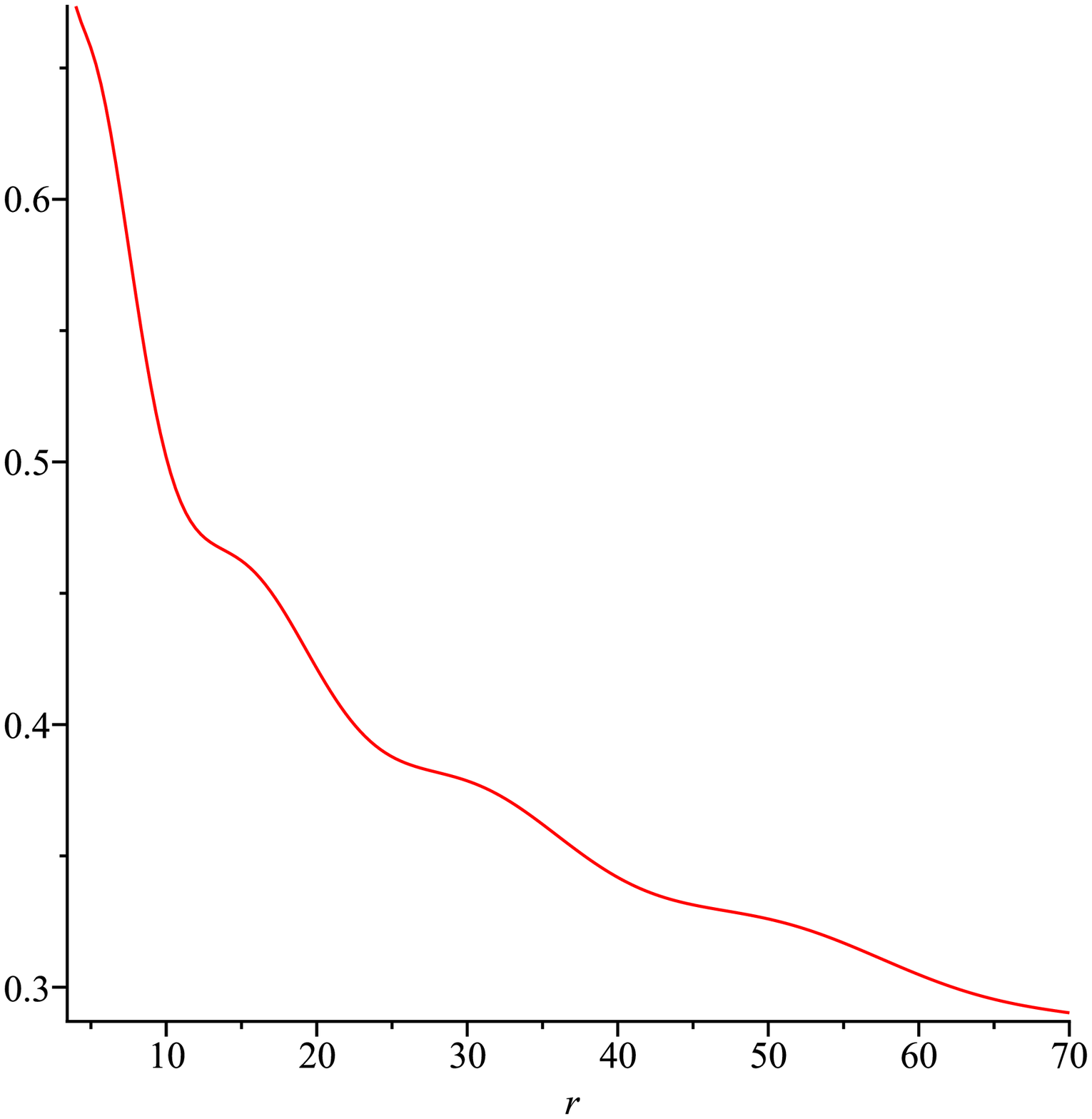 ,width=5cm}} \hspace{0.5cm}{\textbf{Fig.\,6c}\ \ \small Plot of  $B(r)$} \end{minipage}\quad\ \
\begin{minipage}{5cm}
{\psfig{figure=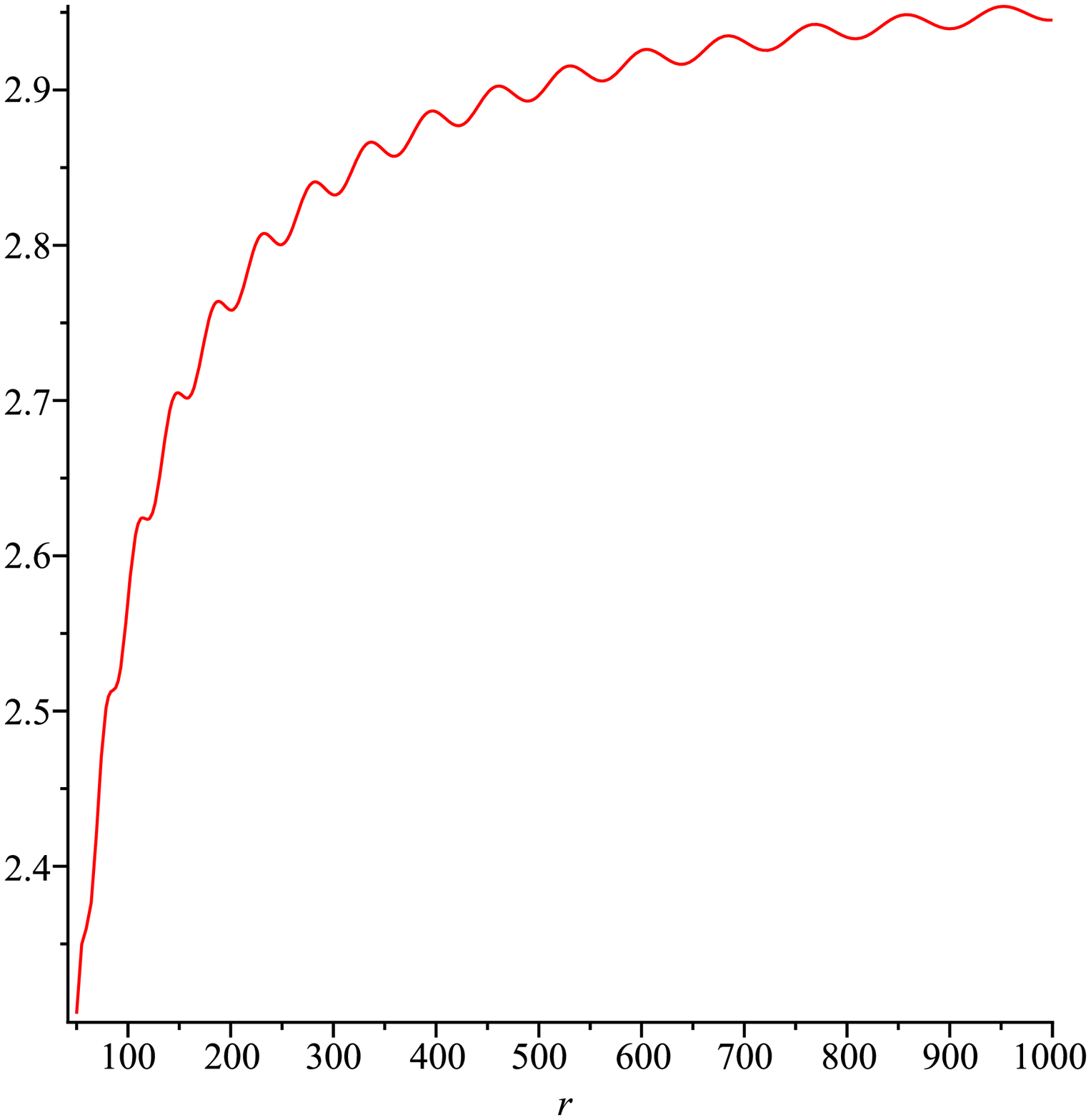,width=5cm}} \hspace{0.5cm}{\textbf{Fig.\,6d}\ \ \small Plot of  $r^{\frac 1 2} B(r)$} \end{minipage}
\end{center}

\begin{center}
\begin{minipage}{5cm}
 {\psfig{figure=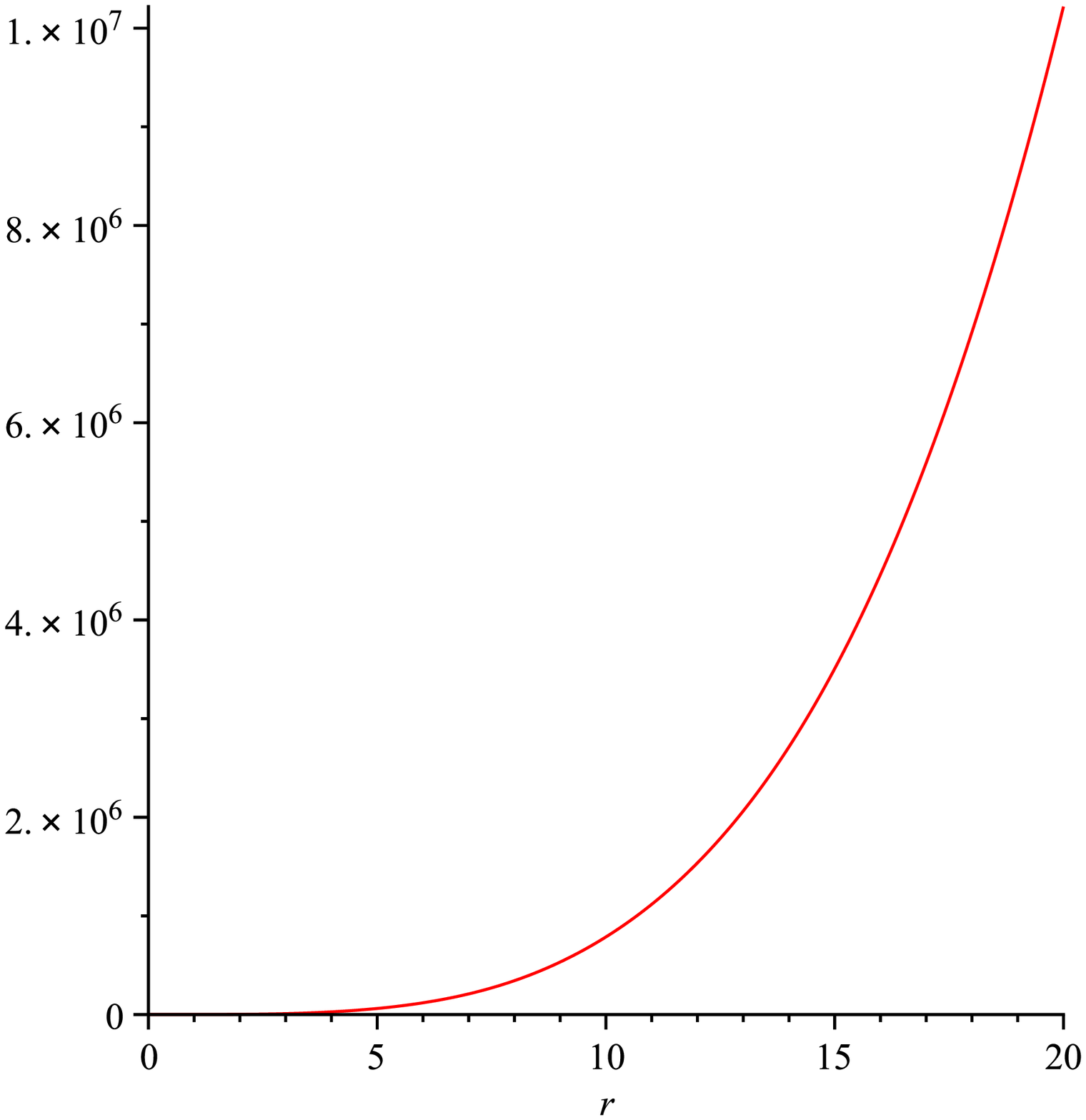 ,width=5cm}} \hspace{0.5cm}{\textbf{Fig.\,6e}\ \  \small  Plot of  $l_3(r)$} \end{minipage}\quad\ \
\begin{minipage}{5cm}
{\psfig{figure=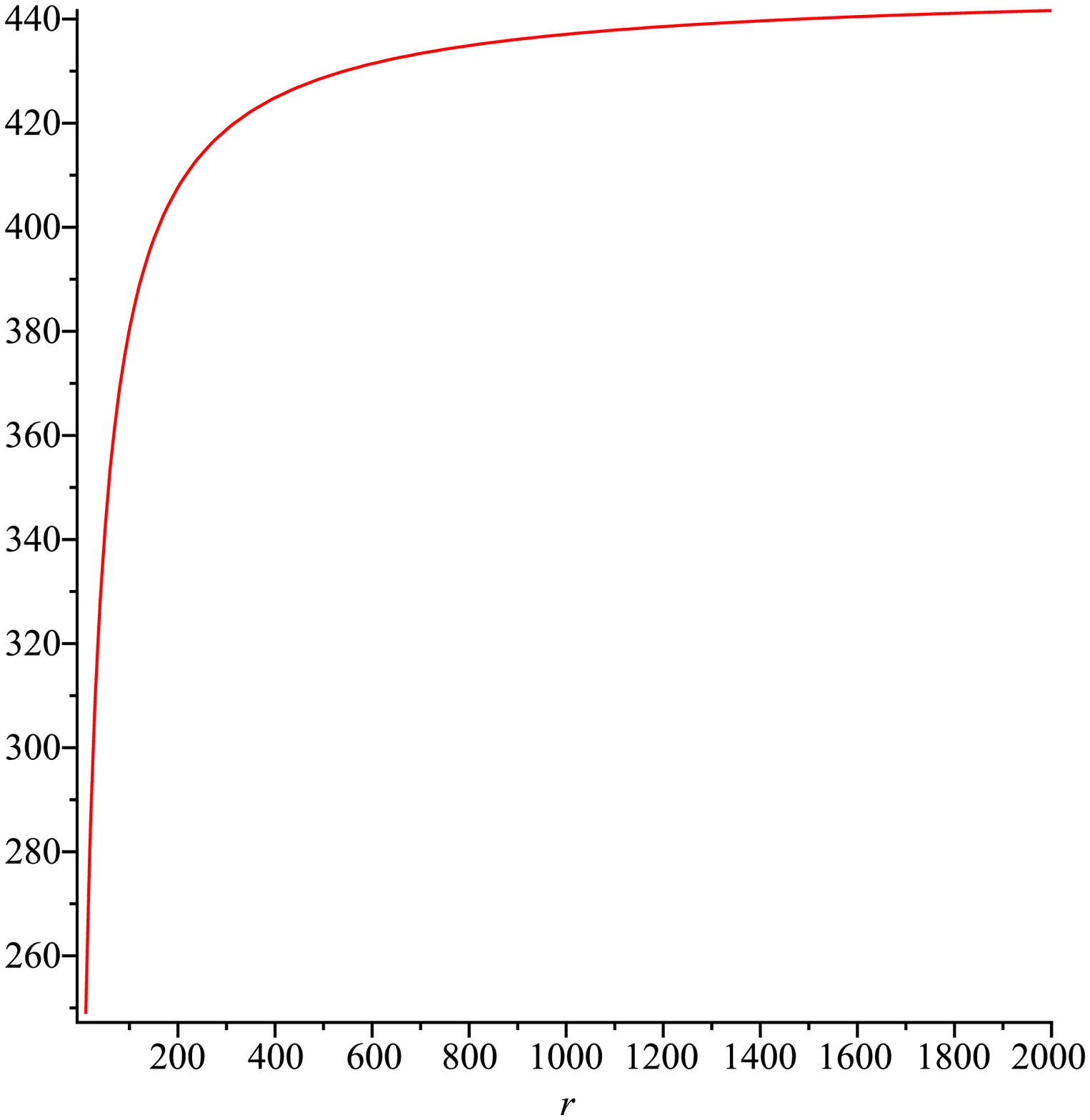,width=5cm}} \hspace{0.5cm}{\textbf{Fig.\,6f}\ \ \small Plot of  $r^{-\frac 7 2} l_3(r)$} \end{minipage}
\end{center}
\end{example}

Example~\ref{ex5} demonstrates that it might be reasonable to use $OR$ class instead of $\cal{L}$ to obtain more general Abelian and Tauberian theorems.  Indeed, such results  were obtained by \citealt{ole1, ole2}.

We will write $f\asymp g$ if $f = O(g)$ and $g = O(f).$ It is convenient to introduce the following notations: \[\tilde{l}_n(r):=\frac{l_n(r)}{r^{2(n-1)}},\quad\quad
\tilde{b}_n(r):=\frac{b_n(r)}{r^{2n}}.\]

\begin{theorem} Let $0>\alpha\ge\beta>2-n.$ Then the following statements are equivalent:
\begin{itemize}
\item[\rm\textbf{(i)}] $G(1/\cdot)\in OR(\beta,\alpha);$
\item[\rm\textbf{(ii)}] $\tilde{l}_n(\cdot)\in OR(\beta,\alpha);$
\item[\rm\textbf{(iii)}] $\tilde{l}_n(r)\asymp G(1/r)$ as $r\to
\infty,$ and there exist positive numbers $C,$ $C'$ and $r_0$ such that
    \[C'\left(\frac{r_1}{r}\right)^\beta\le\frac{G(1/r_1)}{\tilde{l}_n(r)}\le C\left(\frac{r_1}{r}\right)^\alpha,
    \quad r_1\ge r\ge r_0.\]
    \end{itemize}
\end{theorem}

\begin{theorem}\label{th8} Let $0>\alpha\ge\beta>-n.$ Then the following statements are equivalent:
\begin{itemize}
\item[\rm{\textbf{(i)}}] $G(1/\cdot)\in OR(\beta,\alpha);$
\item[\rm\textbf{(ii)}] $\tilde{b}_n(\cdot)\in OR(\beta,\alpha);$
\item[\rm\textbf{(iii)}] $\tilde{b}_n(r)\asymp G(1/r)$ as $r\to
\infty,$ and there exist positive numbers $C,\,C'$ and $r_0$ such
that
\[C'\left(\frac{r_1}{r}\right)^\beta\le\frac{G(1/r_1)}{\tilde{b}_n(r)}\le C\left(\frac{r_1}{r}\right)^\alpha,
    \quad r_1\ge r\ge r_0.\]
    \end{itemize}
\end{theorem}

To deal with spectral densities one can use the following results.

\begin{theorem}\label{th9}  Let $\alpha<n.$ If $g(1/\cdot)\in OR(\beta,\alpha)$ then
$G(1/\cdot)\in
OR(\beta-n,\alpha-n)$ and $G(1/r)\asymp {g(1/r)}/{r^n},$
$r\to \infty.$
\end{theorem}

\begin{definition}  An isotropic spectral density $g(\cdot)$ is called essentially positive (in a neighborhood
of the origin) if there exists a number $\varepsilon >0,$ such that $g(x)>0$ for all $x\in (0;\varepsilon].$
\end{definition}

\begin{theorem} Let $g(1/\cdot)\in OR$ be an essentially positive isotropic spectral density and let
$\alpha<0.$ If $G(1/\cdot)\in OR(\beta,\alpha),$ then $g(1/\cdot)\in OR(\beta+n,\alpha+n)$ and
$G(1/r)\asymp {g(1/r)}/{r^n},$ as $r\to \infty.$
\end{theorem}

Abelian and Tauberian theorems for $O$-regularly varying spectral densities were obtained by \citealt{ole1}. $OR$ relations of the local behaviour of the spectral distribution function $G(\cdot)$ in a neighborhood of an arbitrary point  $a\in [0,+\infty)$ and the variance of some weighted integral transforms of random fields were investigated by \citealt{ole2}.

The last example in this section shows that Abelian and Tauberian theorems for $OR$ class can not be reduced to narrower $\cal{L}$ cases considered before. We give an example of functions $g(1/\cdot)$ and corresponding $\tilde{b}_n(\cdot)$ which both belong to $OR$ class, but are not in $\cal{L}.$
\begin{example}
Let  the spectral distribution function be determined by its derivative

\begin{equation} \label{prOR}     G'(\lambda)= \left\{%
\begin{array}{ll}
        0 , &   \lambda> \frac{1}{{T_n}};\\
        1,  &     \lambda\in\left. \left(\frac{1}{2T_n^{2k+1}},  \frac{1}{T_n^{2k+1}}  \right.  \right], \; k\ge 0\, ;\\
        \varepsilon,  &     \lambda\in\left. \left(\frac{1}{T_n^{2(k+1)+1}},  \frac{1}{2T_n^{2k+1}}    \right. \right],
      \;\;k\ge 0\,  ,\\
\end{array}%
\right.
\end{equation}
where  $T_n>2$ such that
\begin{equation} \label{integralSn}  \int_{\frac{T_n}{2}}^\infty S_n(\lambda)d\lambda<\delta_1,\;\;    \int_{0}^\frac {1}{T_n}S_n(\lambda)d\lambda<\delta_2,
\end{equation}
and $S_n(\lambda):={J_{\frac n2}^2(\lambda)}/{\lambda^n}.$ $S_n(\lambda)$ is an integrable function. Therefore, for any positive $\delta_1$ and $\delta_2$ there exists $T_n$ such that \eqref{integralSn} holds.

Positive constants $\varepsilon,$ $\delta_1,$ and $\delta_2$ will be defined later.

Let $u:= 1/\lambda, \;\; \tilde{G}'(u):= G'(1/u)=G'(\lambda).$
By \eqref{prOR} for $t>1$ and $u\ge 1$ we get
$$\varepsilon  \le \frac{\tilde{G}'(t u)}{\tilde{G}'(u)}\le \frac 1 \varepsilon. $$
Hence, $G'(1/\cdot)\in OR(0,0).$ However, $G'(1/\cdot)\not\in \cal{L}$ because  $\lim_{u\to\infty}{\tilde{G}'(t u)}/{\tilde{G}'(u)}$ does not exist for $t>1.$ By Definition~\ref{def3} we obtain that $g\left(1/\cdot\right)\in OR(n-1, n-1),$ but $g\left(1/u\right)/u^{n-1}\not\in \cal{L}.$

Note that $\tilde{b}_n(r)\in OR(-1,-1)$ by Theorems~\ref{th8} and \ref{th9}.

Let us show that $r\tilde{b}_n(r)\not \in \cal{L}.$ The substitution of the sequence $\lambda_k=T_n^{2k},$ $k\ge 1$ and $t=T_n$ in Definition~\ref{def4} gives

$$   \frac{\tilde{b}_n(\lambda_k t)}{\tilde{b}_n(\lambda_k)}=\frac 1 t \cdot
\frac{\int_0^\infty S_n(\lambda)G'\left(\frac{\lambda}{\lambda_k t}\right)d\lambda }{\int_0^\infty S_n(\lambda)G'\left(\frac{\lambda}{\lambda_k }\right)d\lambda} \ge     \frac 1 t \cdot  \frac{A}{\delta_1+\delta_2+B\cdot \varepsilon },  $$
where
$$A:= \int_{\frac 1 2}^1 S_n(\lambda)d\lambda, \;\; B:=\int_0^\infty S_n(\lambda) d\lambda.$$
Note that both $A$ and $B$ are finite positive numbers.

If $\lambda_k=T_n^{2k+1},$ $k\ge 0,$ and $t=T_n,$ then
$$   \frac{\tilde{b}_n(\lambda_k t)}{\tilde{b}_n(\lambda_k)}=\frac 1 t \cdot
\frac{\int_0^\infty S_n(\lambda)G'\left(\frac{\lambda}{\lambda_k t}\right)d\lambda }{\int_0^\infty S_n(\lambda)G'\left(\frac{\lambda}{\lambda_k }\right)d\lambda} \le     \frac 1 t \cdot  \frac{\delta_1+\delta_2+B\cdot \varepsilon }{A}.  $$

If we choose $\delta_1,$ $\delta_2$ and $\varepsilon$ so small that
$$ \frac {\delta_1+\delta_2+B\cdot \varepsilon }{A}<\frac{A} {\delta_1+\delta_2+B\cdot \varepsilon }$$
then the limit
 $$  \lim_{r\to\infty} \frac{\int_0^\infty S_n(\lambda)G'\left(\frac{\lambda}{r t}\right)d\lambda }{\int_0^\infty S_n(\lambda)G'\left(\frac{\lambda}{r}\right)d\lambda} $$
does not exist. Therefore $r\tilde{b}_n(r)\not \in \cal{L}.$
\end{example}

\section{Radial directional asymptotic behaviour}

In this section we discuss some generalizations of Tauberian and Abelian theorems for "radial directional" homogeneous random fields. Even if a random field is not isotropic, its covariance functions can exhibit some regular asymptotic behaviour in each radial direction, see \citealt{dou} and a proper definition in Theorem~\ref{th11}. This issue has no parallel in  one-dimensional time series modeling.

Limit theorems and statistical applications of such LRD random fields defined on the discrete parameter space  $\mathbb{Z}^{n}$ were discussed by \citealt{dou}, \citealt{lav1,lav2,lav3,lav4}, \citealt{ber}, etc. Some examples of "directional" homogeneous LRD random fields on $\mathbb{Z}^{n}$ were given.

To show the presence of a LRD behaviour mentioned papers used Tauberian and Abe\-lian theorems for Fourier transforms by \citealt{wai}. We adjust Wainger's results to LRD random fields defined on the continuous parameter space~$\mathbb{R}^{n}.$

\begin{definition} An infinitely differentiable function $L(\cdot)$ belongs to the class $\mathcal{\tilde{L}}$ if
\begin{enumerate}
  \item for any $\delta>0$ there exist  $\lambda_0(\delta)>0$ such that $\lambda^{-\delta}L(\lambda)$ is decreasing and $\lambda^{\delta}L(\lambda)$ is increasing if $\lambda>\lambda_0(\delta);$
  \item $L_j(\cdot)\in \mathcal{L}$ for all $j\ge 0,$ where $L_0(\lambda):=L(\lambda),$ $L_{j+1}(\lambda):=\lambda L_{j}'(\lambda).$
\end{enumerate}
\end{definition}

 Any function $S(\cdot)\in C^\infty(s_{n-1}(1))$ possesses the mean-square representation
$$S(\omega)=\sum_{k=0}^\infty\sum_{m=-k}^k a(k,m)Y_{k,m}(\omega),\quad \omega\in s_{n-1}(1),$$
where $Y_{k,m}(\cdot)$ is a spherical harmonic of degree $k,$ see \citealt{andr}.

The function
$$\tilde{S}_{\alpha,n}(\omega):=\pi^{\alpha-n}\sum_{k=0}^\infty\sum_{m=-k}^k(-i)^k\, a(k,m)\frac{\Gamma\left(\frac{n+k-\alpha}{2}\right)}{\Gamma\left(\frac{k+\alpha}{2}\right)}Y_{k,m}(\omega)$$
is well defined in the mean-square sense and belongs to $C^\infty(s_{n-1}(1)).$

\begin{theorem}\label{th11} Let $\alpha\in (0,n),$ $S(\cdot)\in C^\infty(s_{n-1}(1)),$ and $L(\cdot)\in \mathcal{\tilde{L}}.$  Let $\xi (x),$ $x\in \mathbb{R}^{n},$ be a mean-square continuous homogeneous random field with mean zero. Let the field $\xi (x)$ have the spectral density $f(u),$ $u\in R^{n},$ which is infinitely differentiable for all $u\not = 0.$

If the covariance function $B(x),$ $x\in R^{n},$ of the field has the following behaviour
\begin{enumerate}
  \item[\rm{(a)}] $ ||x||^\alpha B(x)\sim  S\left(\frac{x}{||x||}\right) L\left(||x||\right),$ $||x||\to +\infty,$
\end{enumerate}
then the spectral density satisfies the condition
 \begin{enumerate}
  \item[\rm{(b)}] $||u||^{n-\alpha} f(u)\sim  \tilde{S}_{\alpha,n}\left(\frac{u}{||u||}\right) L\left(\frac{1}{||u||}\right),$ $||u||\to 0.$
     \end{enumerate}
\end{theorem}

We begin with an example of a spectral density and a covariance function which illustrates Theorem~\ref{th11}.
\begin{example}\label{ex7}
Let $n=3,$ $\alpha=2,$ and $\xi(x)$ be a mean-square continuous homogeneous random field with the spectral density (compare with Example~\ref{ex2})
\begin{equation}\label{fu}
f(\rho, \beta, \psi)=(4+3\cos^2(\beta))\cdot\frac{K_{\frac{1}{2}}(\rho)}{\sqrt{2^3\pi^3\,\rho}}\,,
\end{equation}
where $(\rho, \beta, \psi),$ $\rho\ge 0,$  $\beta\in [0,\pi],$ $\psi \in [0,2\pi),$ are spherical coordinates of the point~$u\in\mathbb{R}^3.$

 Using the spherical harmonics $Y_{0,0}\left(\frac{u}{||u||}\right):=1$ and $Y_{2,0}\left(\frac{u}{||u||}\right):=\frac{3\cos^2(\beta)-1}{2}$
we can represent (\ref{fu}) as
\begin{equation}\label{fu0}
f(u)=\left(5Y_{0,0}\left(\frac{u}{||u||}\right)+2Y_{2,0}\left(\frac{u}{||u||}\right)\right)\cdot\frac{K_{\frac{1}{2}}(||u||)}{\sqrt{2^3\pi^3\,||u||}}.
\end{equation}

Hence, the spectral decomposition~(\ref{Bx}) can be rewritten as
$$
B(x) =4\pi\int_0^\infty\frac{\rho^2 K_{\frac{1}{2}}(\rho)}{\sqrt{2^3\pi^3\,\rho}}\int_{s_{2}(1)}e^{i||x||\rho\left\langle \frac{x}{||x||},\omega  \right\rangle}\left(5Y_{0,0}\left(\omega\right)+2Y_{2,0}\left(\omega\right)\right)
d\sigma(\omega)\ d\rho\,.
$$

Using properties of the Fourier transform of spherical harmonics, see \citealt[(9.10.2)]{andr}, one may evaluate the internal surface integral over the sphere $s_{2}(1)$ in the following way
$$
B(x) =4\pi\left(5Y_{0,0}\left(\frac{x}{||x||}\right)\int_0^\infty\frac{J_{1/2}(\rho||x||)}{ \sqrt{\rho||x||}}\rho^{3/2}K_{\frac{1}{2}}(\rho)d\rho\right.$$
$$\left.-2Y_{2,0}\left(\frac{x}{||x||}\right)\int_0^\infty\frac{J_{5/2}(\rho||x||)}{ \sqrt{\rho||x||}}\rho^{3/2}K_{\frac{1}{2}}(\rho)d\rho\right)
\,.
$$
Finally we obtain
$$
B(r, \theta, \varphi) =\frac{4\pi}{\left( 1+{r}^{2}\right){r}^{3}} \, \left( {r}^{3}(7-6\, \cos^{2}(\theta))+3{r}^{2}(3\,  \cos^{2}(\theta) -1) \arctan (r) \right.$$
\begin{equation}\label{dir1}
\left.+3\,r(1-3\, \cos^{2} ( \theta )) +3(3 \,  \cos^{2} (\theta) -1) \arctan(r)  \right)\,,\end{equation}
where $(r, \theta, \varphi),$ $r\ge 0,$  $\theta\in [0,\pi],$ $\varphi \in [0,2\pi),$ are spherical coordinates of the point~$x\in\mathbb{R}^3.$

It is obvious from (\ref{dir1}) that the statement (a) of Theorem~\ref{th11} holds, i.e.
$$ ||x||^2 B(x)\sim  S\left(\frac{x}{||x||}\right) L\left(||x||\right) ,\quad ||x||\to +\infty,$$
where $S(\omega)=5 Y_{0,0}(\omega)-4\,\pi\, Y_{2,0}(\omega)=7-6\, \cos^{2}(\theta)$ and $L(\cdot)\equiv const.$

Therefore $\tilde{S}_{\alpha,n}(\omega)=\pi^{-1/2}\left(5 Y_{0,0}(\omega)+2 Y_{2,0}(\omega)\right)$ and the density defined by equation (\ref{fu0}) satisfies the statement (b) of Theorem~\ref{th11}.

The correlation function $B(x)$ does not depend on $\varphi \in [0,2\pi).$ Therefore, we plot it as a function of the variables $r$ and $\theta.$ To apparently display the function $B(x)$ we extend it to an even function on the interval $\theta\in (-\pi,\pi).$  Fig.~7 presents plots of the covariance function $B(x)$  and its normalized transformation in the cylindrical coordinate system $(r\cos(\theta), r\sin(\theta), B(r, \theta, \varphi)).$\newpage
\noindent\begin{minipage}{5.8cm}
 {\psfig{figure=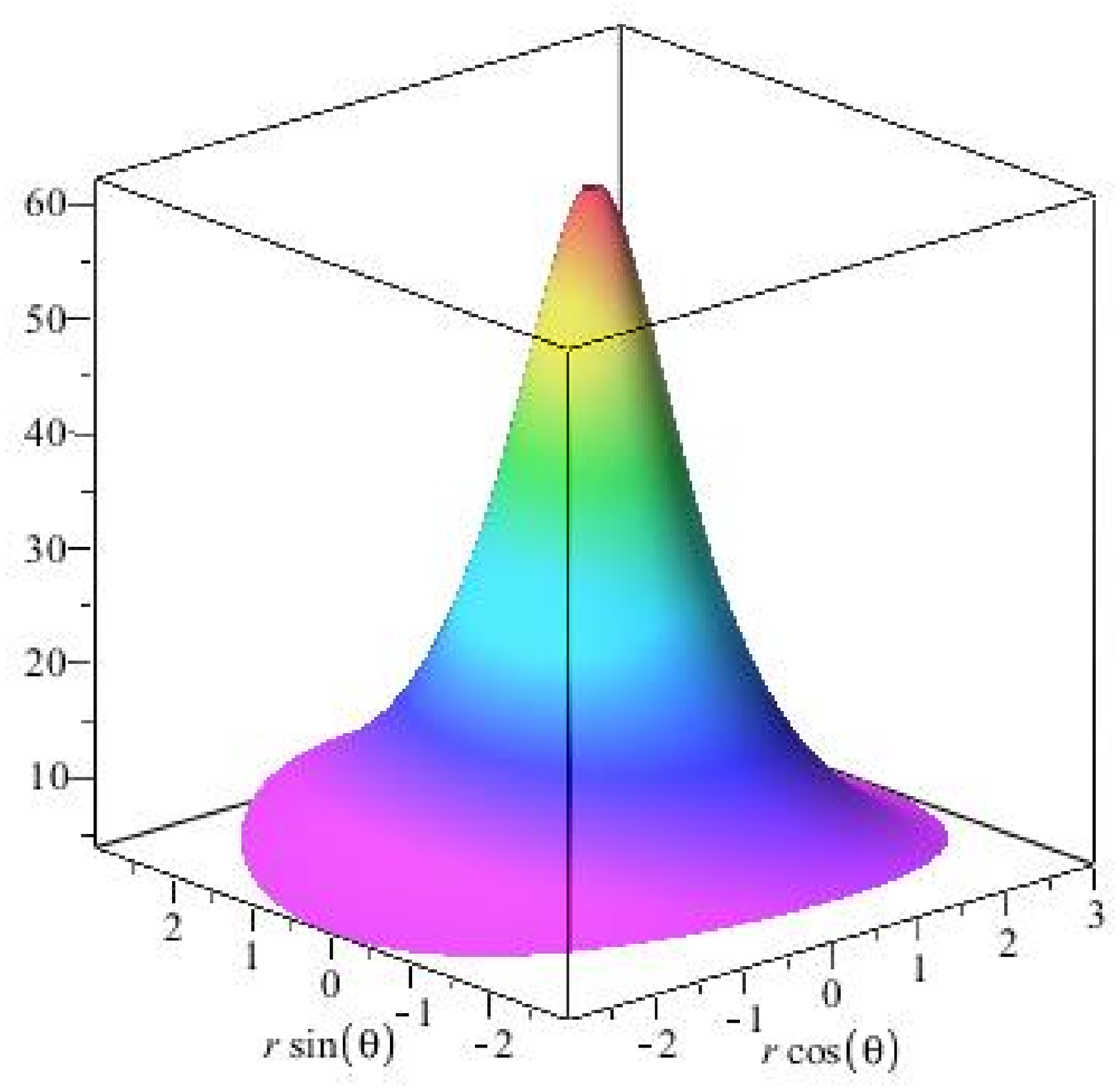,width=6cm}} \vspace{-3mm}\hspace{1cm}{\textbf{Fig.\,7a}\ \
 \small Plot of $B(r, \theta, \varphi)$} \end{minipage}\quad\ \
\begin{minipage}{5.8cm}
{\psfig{figure=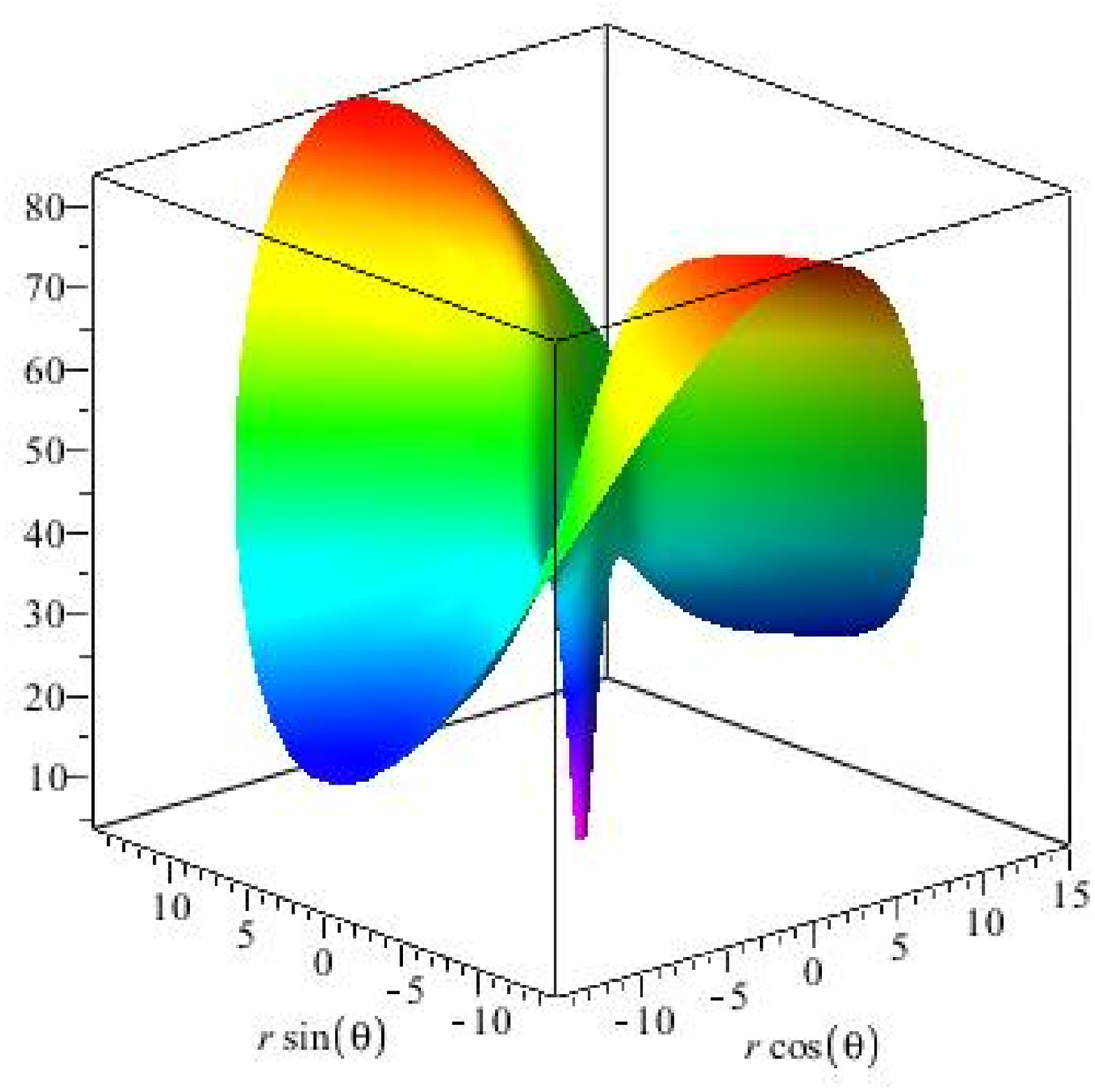,width=6cm}} \vspace{-3mm}\hspace{1cm}{\textbf{Fig.\,7b}\ \  \small Plot of
$r^2B(r, \theta, \varphi)$} \end{minipage}

\noindent\begin{minipage}{5.8cm}
 {\psfig{figure=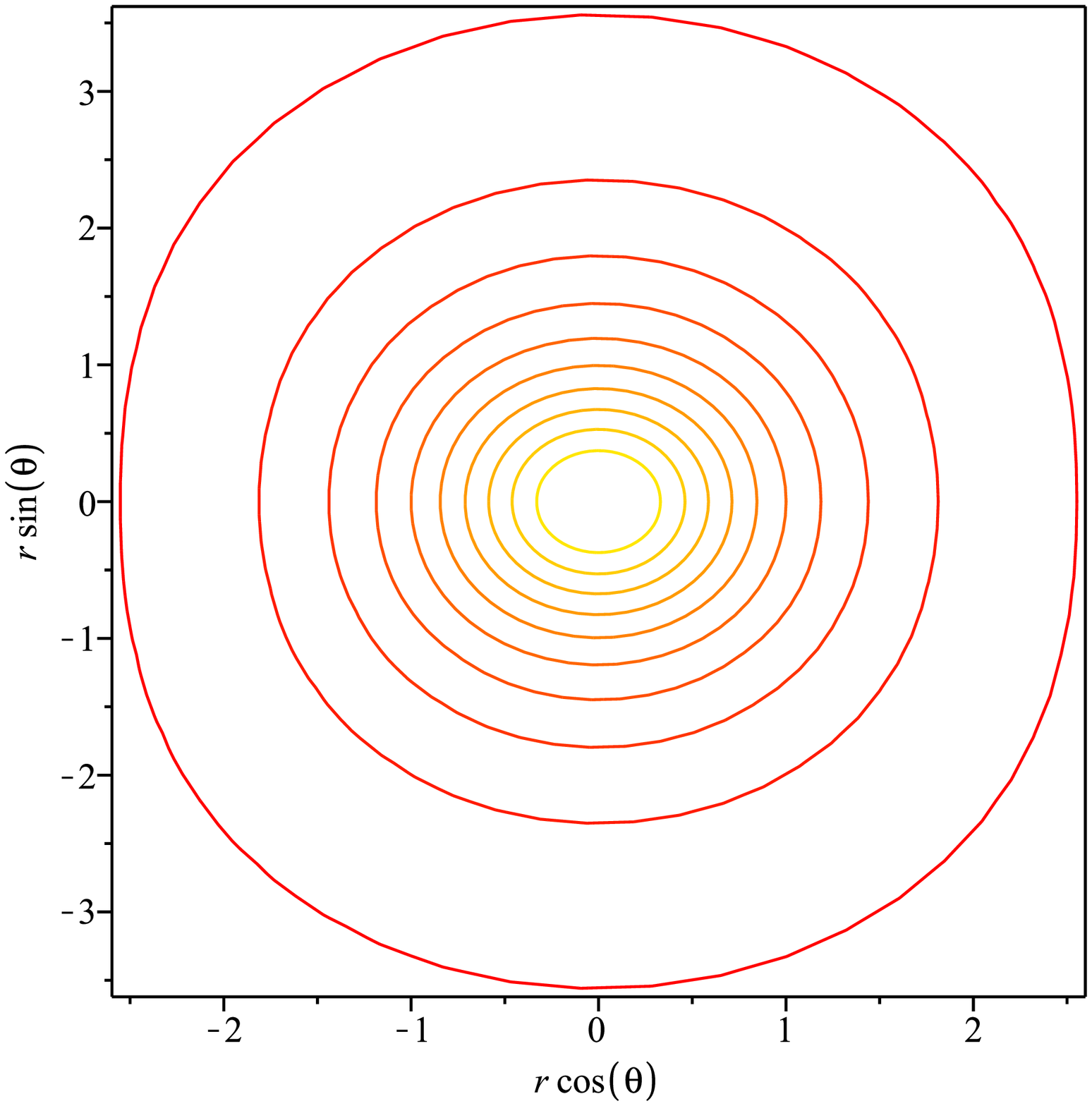,width=5.8cm}} \hspace{0.5cm}{\textbf{Fig.\,8a}\ \
 \small Contour plot of $B(r, \theta, \varphi)$ for\\${}^{}$\hspace{1.4cm} the  cylindrical coordinate system} \end{minipage}\quad\ \
\begin{minipage}{5.8cm}
{\psfig{figure=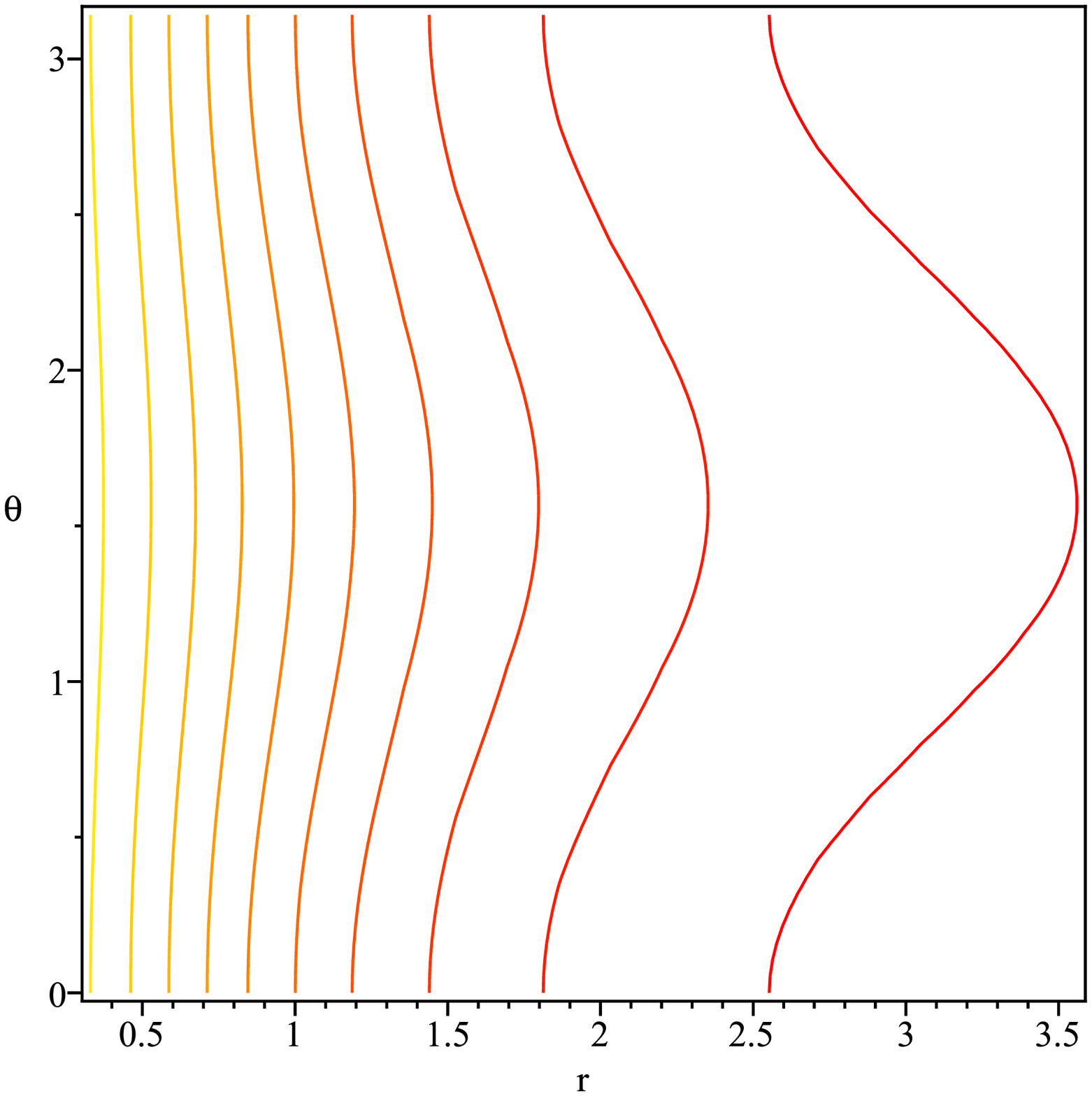,width=5.8cm}} \hspace{0.5cm}{\textbf{Fig.\,8b}\ \  \small Contour plot of $B(r, \theta, \varphi)$ in\\${}^{}$\hspace{1.7cm} $(r, \theta)$  coordinates} \end{minipage}

{}\vspace{-2mm}\noindent\begin{minipage}{5.8cm}
 {\psfig{figure=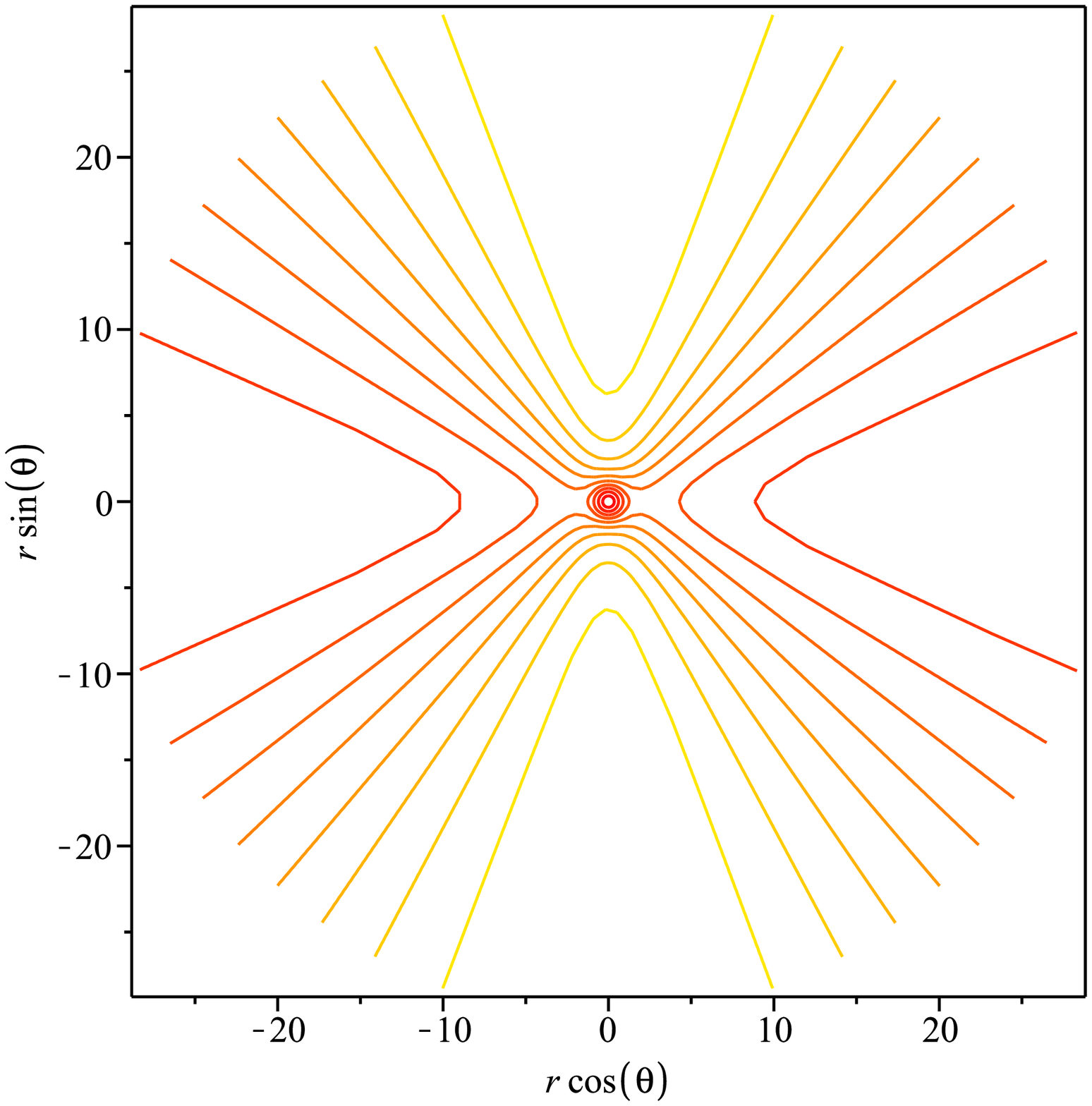,width=5.8cm}} \hspace{0.5cm}{\textbf{Fig.8c}\
 \small Contour plot of $r^2B(r, \theta, \varphi)$ for\\${}^{}$\hspace{1.4cm}  the cylindrical coordinate system} \end{minipage}
 \quad\ \
\begin{minipage}{5.8cm}
{\psfig{figure=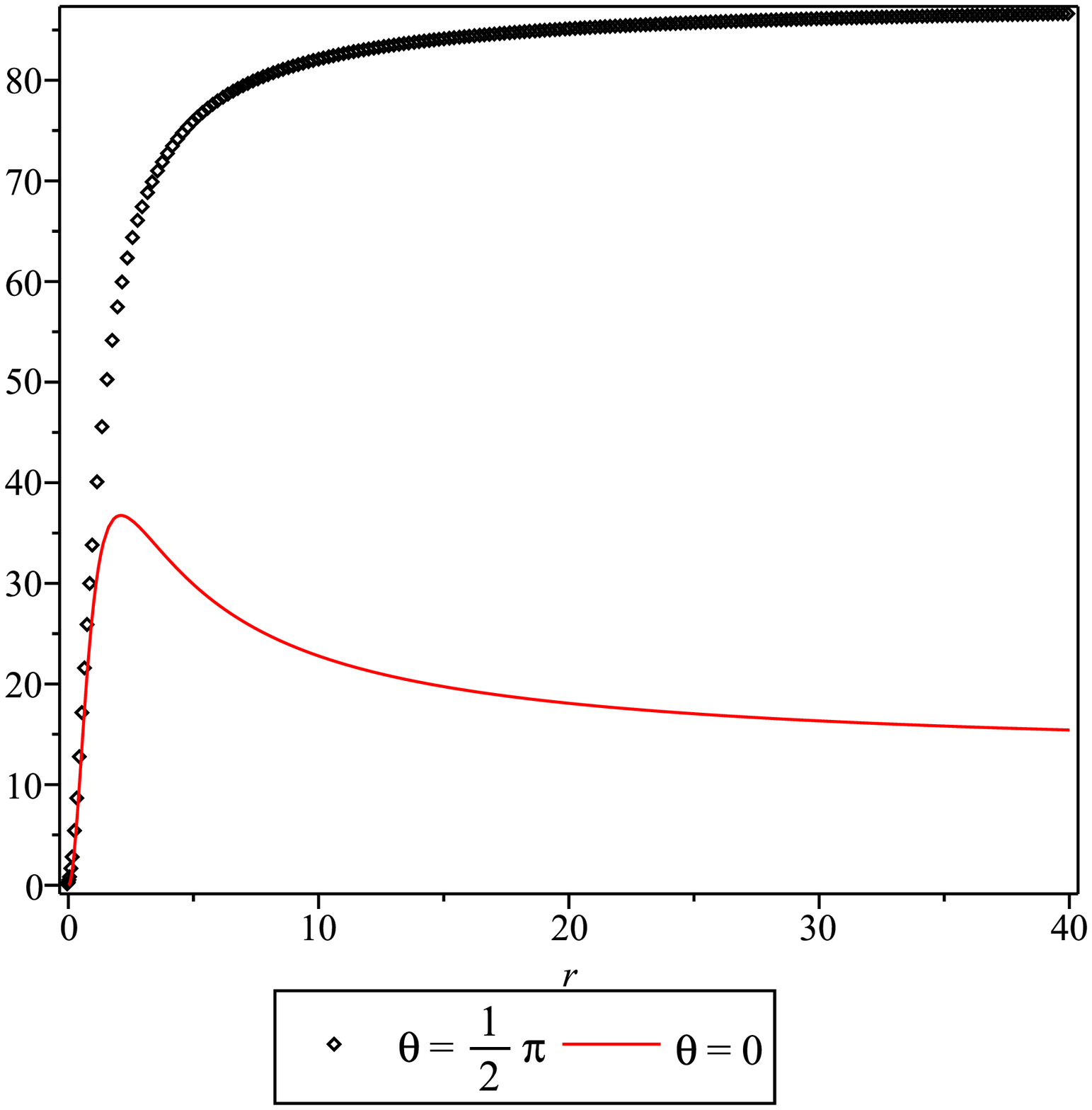,width=5.8cm}} \hspace{0.5cm}{\textbf{Fig.\,8d}\ \  \small Plots of $r^2B(r, \theta, \varphi)$ for \\${}^{}$\hspace{1.7cm}  $\theta=0$ and $\theta=\pi/2.$} \end{minipage}

The contour plots in Fig.~8a-8c demonstrate that random field $\xi(x)$ exhibits "directional" homogeneous behaviour.
\end{example}

\begin{definition} A random field $\xi (x),$ $x\in \mathbb{R}^{n},$ is anisotropic if its covariance function $B(x)=B_0(||Ax||)$ for some isotropic covariance function $B_0(\cdot)$ and nondegenerate matrix $A.$
\end{definition}

\begin{remark}
It seems that the contour plots in Fig.8a and Fig.8b demonstrate anisotropic behaviour of $\xi(x).$  However, it is not true in the general case, because random fields with covariance functions satisfying the condition~(a) of Theorem~\ref{th11} are "radially directional" homogeneous, but not anisotropic. Even more, the example below shows the random fields with a radial spectral distribution function satisfying the condition~(b) of Theorem~\ref{th11}. In spite of that, the field is neither "radially directionally" homogeneous nor anisotropic in terms of covariance functions.
\end{remark}

The following example shows that the conditions of Theorem~\ref{th11} are not necessary and, similarly to Example~\ref{ex5},  the statements (a) and (b) in Theorem~\ref{th11} are not equivalent in the general case.

\begin{example} We continue with the notations and transformations of Example~\ref{ex7} and consider a new spectral density.

Let $n=3,$ $\alpha=2,$ and the spectral density is
$$
f(\rho, \beta, \psi)=
\left\{
  \begin{array}{ll}
    \frac{4+3\cos^2(\beta)}{2^{3/2}\pi^{3/2}\rho}\,, & \hbox{if}\ \rho\in (0,1]; \\
    0, & \hbox{otherwise.}
  \end{array}
\right.
$$
Then, the statement (b) of Theorem~\ref{th11} holds.

The corresponding covariance function is
$$
B(r, \theta, \varphi) =\frac {2^{5/2}\sqrt {\pi }}{{r}^{3}} \left( r(7-4\,\cos \left( r
 \right)-3\, \cos^{2} \left( \theta \right)
 (2+ \cos
 \left( r \right)))\right.$$
 $$\left.+3(3\,  \cos^{2} \left( \theta \right)  -1)\sin
 \left( r \right)  \right).
$$
It is obvious that this function does not satisfy the statement (a) of Theorem~\ref{th11}.

\

\noindent\begin{minipage}{5.8cm}
 {\psfig{figure=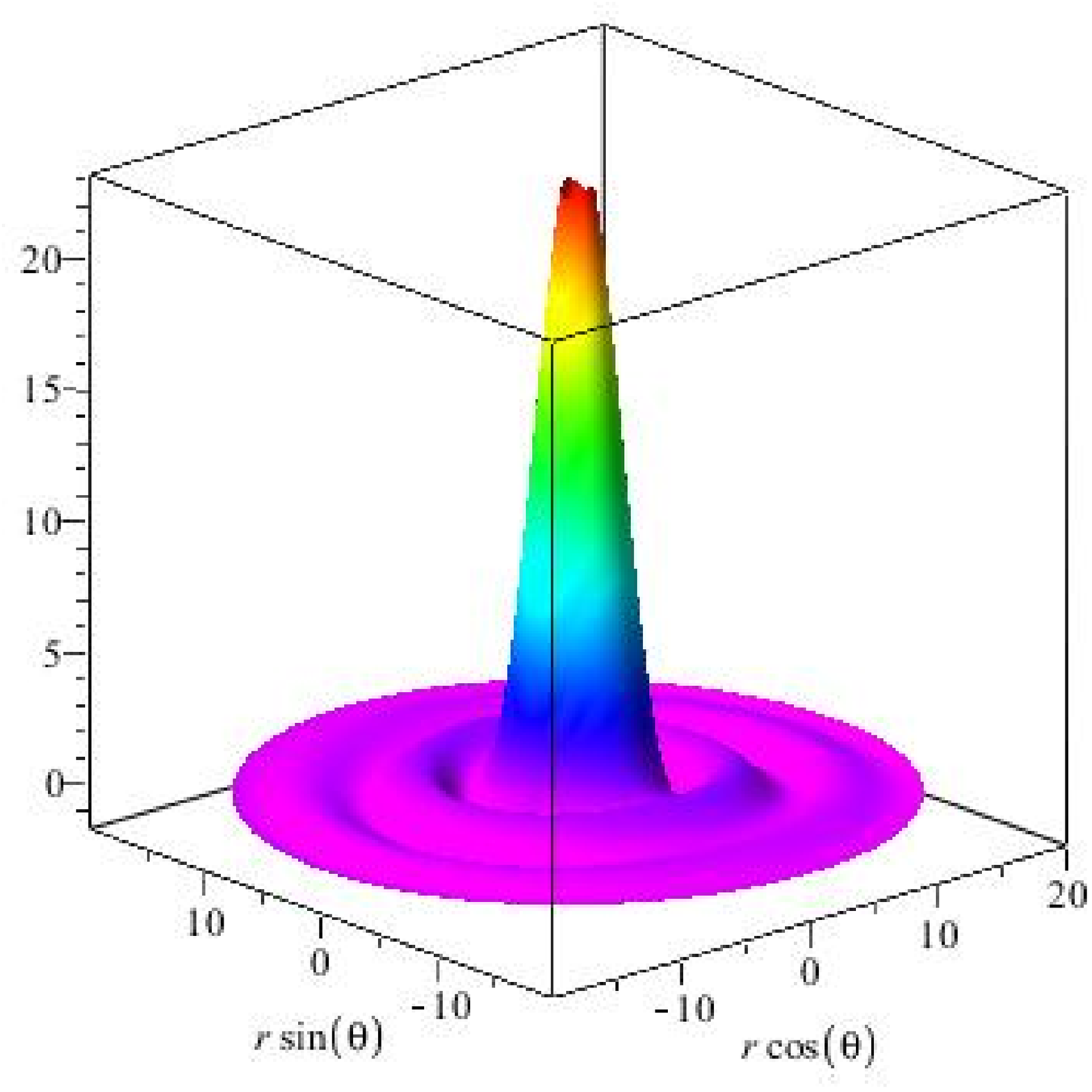,width=5.8cm}} \hspace{1cm}{\textbf{Fig.\,9a}\ \
 \small Plot of $B(r, \theta, \varphi)$} \end{minipage}\quad\ \
\begin{minipage}{5.8cm}
{\psfig{figure=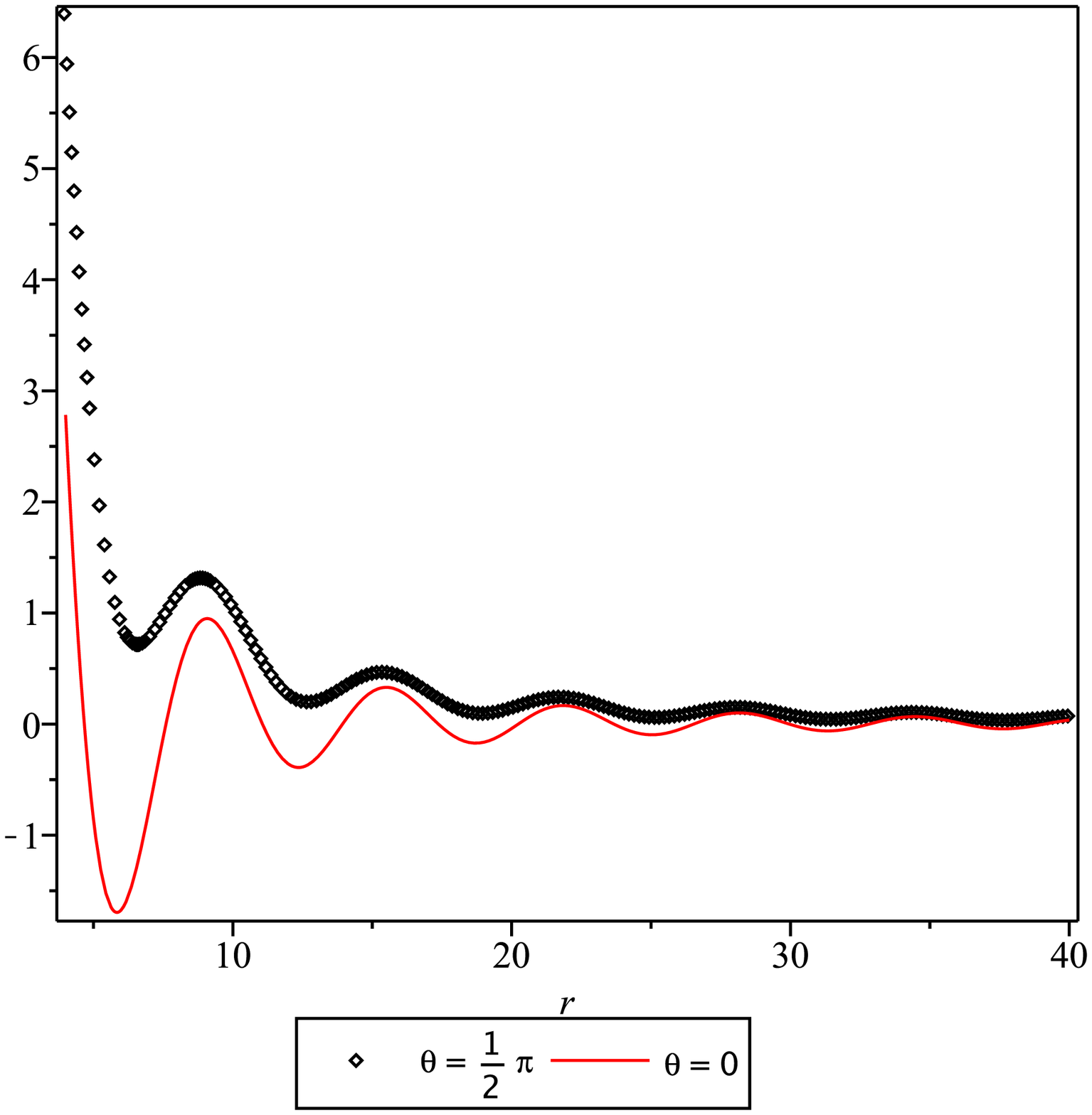,width=5.8cm}} \hspace{1cm}{\textbf{Fig.\,9b}\ \  \small Plots of $B(r, \theta, \varphi)$ for \\${}^{}$\hspace{1.7cm}  $\theta=0$ and $\theta=\pi/2.$} \end{minipage}\vspace{5mm}

\noindent\begin{minipage}{5.8cm}{\psfig{figure=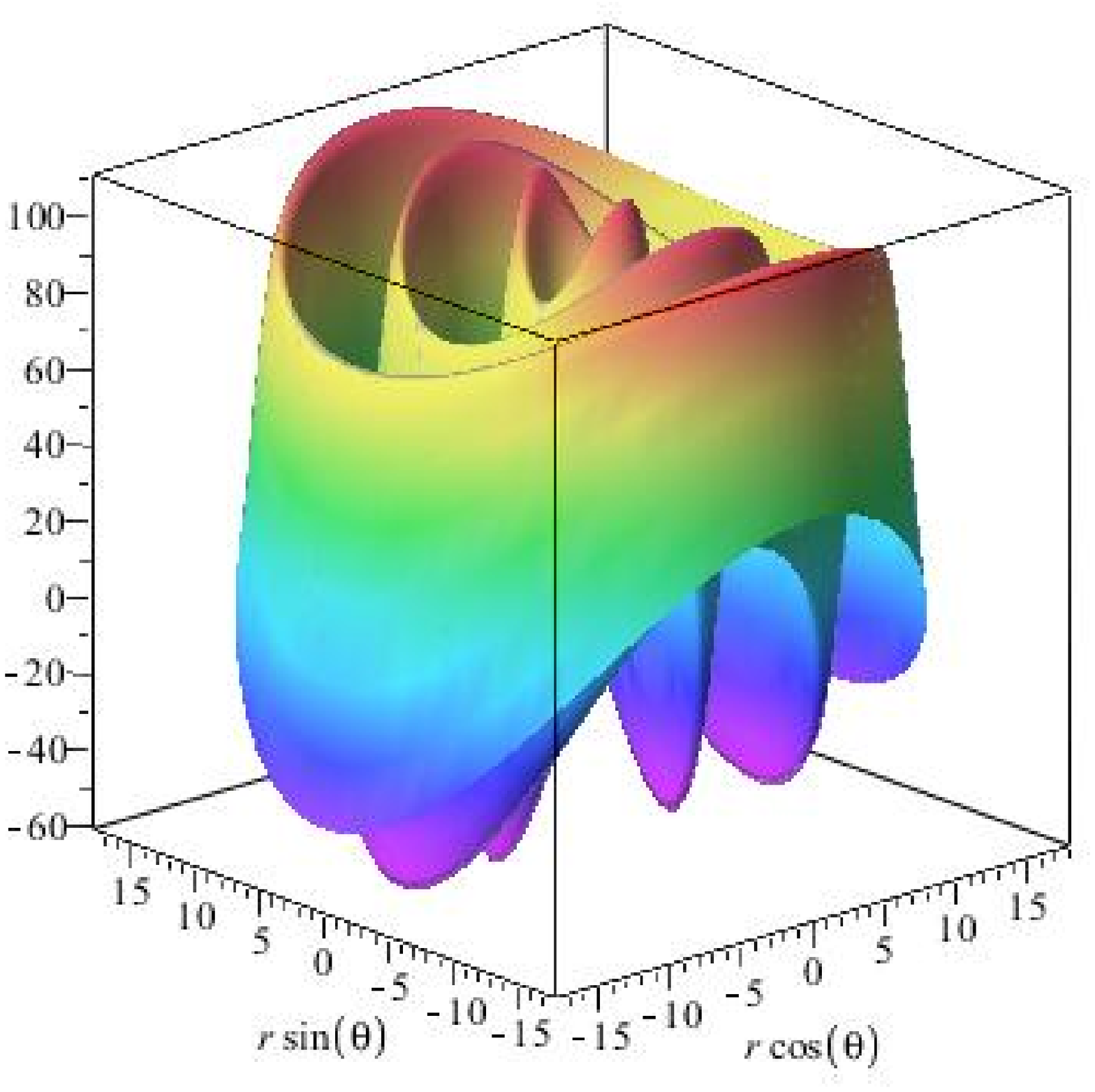,width=5.8cm}} \hspace{1cm}{\textbf{Fig.\,9c}\ \  \small Plot of $r^2B(r, \theta, \varphi)$}\end{minipage}
\quad\ \
\begin{minipage}{5.8cm}
{\psfig{figure=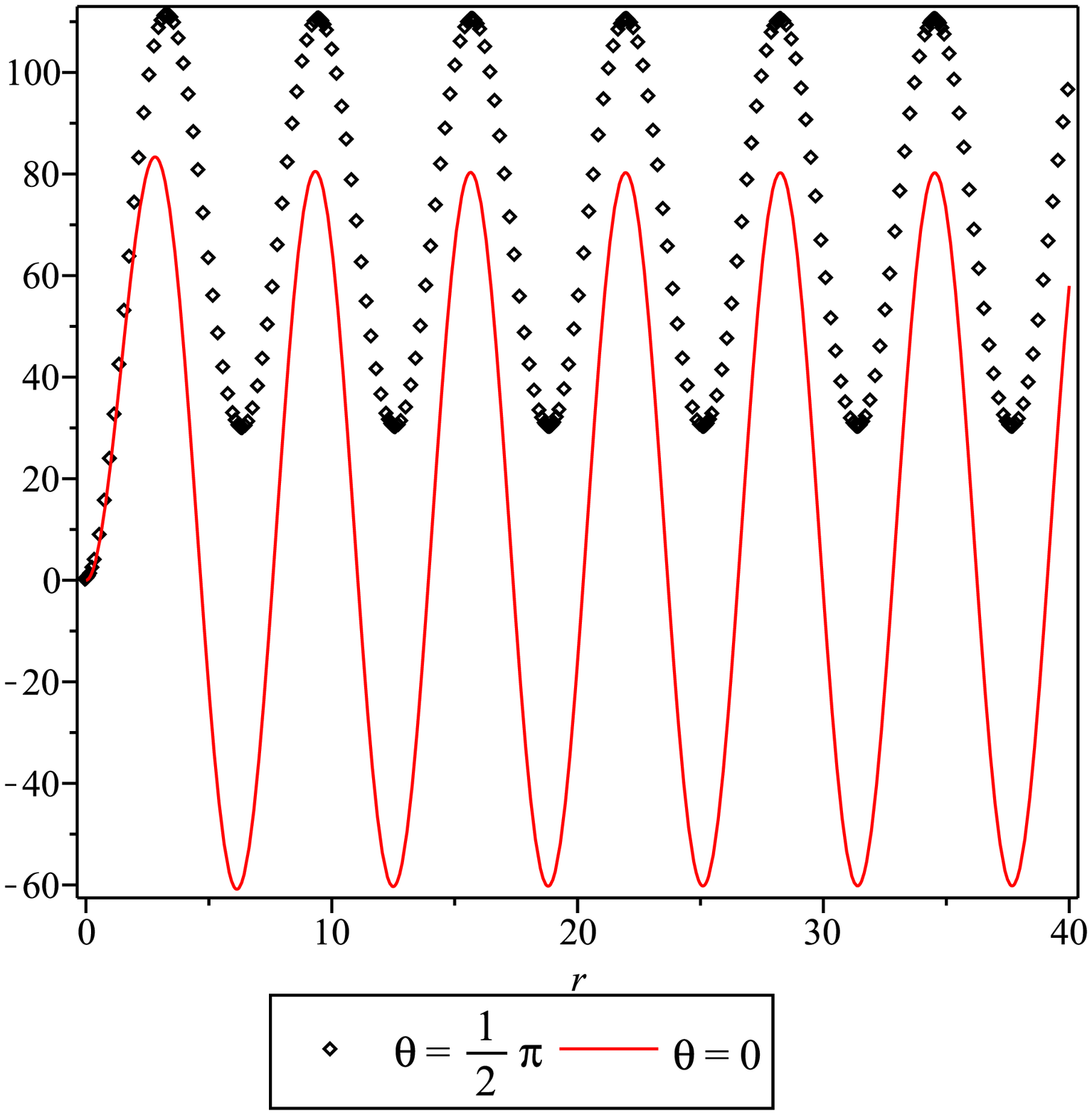,width=5.8cm}} \hspace{0.5cm}{\textbf{Fig.\,9d}\ \  \small Plots of $r^2B(r, \theta, \varphi)$ for \\${}^{}$\hspace{1.7cm}  $\theta=0$ and $\theta=\pi/2.$} \end{minipage}\vspace{5mm}

As in Example~\ref{ex7} we extend  $B(x)$ to an even function on the interval $\theta\in (-\pi,\pi)$ and plot it as a function of the variables $r$ and $\theta.$ Plots of the covariance function $B(x)$  and its normalized transformation in the cylindrical coordinate system $(r\cos(\theta), r\sin(\theta), B(r, \theta, \varphi))$ are shown in Fig.~9.

 It is well-known that anisotropic random fields should have the same covariance structure with different covariance scales along different directions. Directionally homogeneous fields possess the same radial covariance function up to a multiplier, which depends only on the direction. It seems that the contour plots in Fig.~10a and Fig.~10b demonstrate  either anisotropy or "directional" homogeneity of $\xi(x).$ However, two plots in  Fig.~9d and Fig.~10c prove that neither is true.

\noindent\begin{minipage}{5.8cm}
 {\psfig{figure=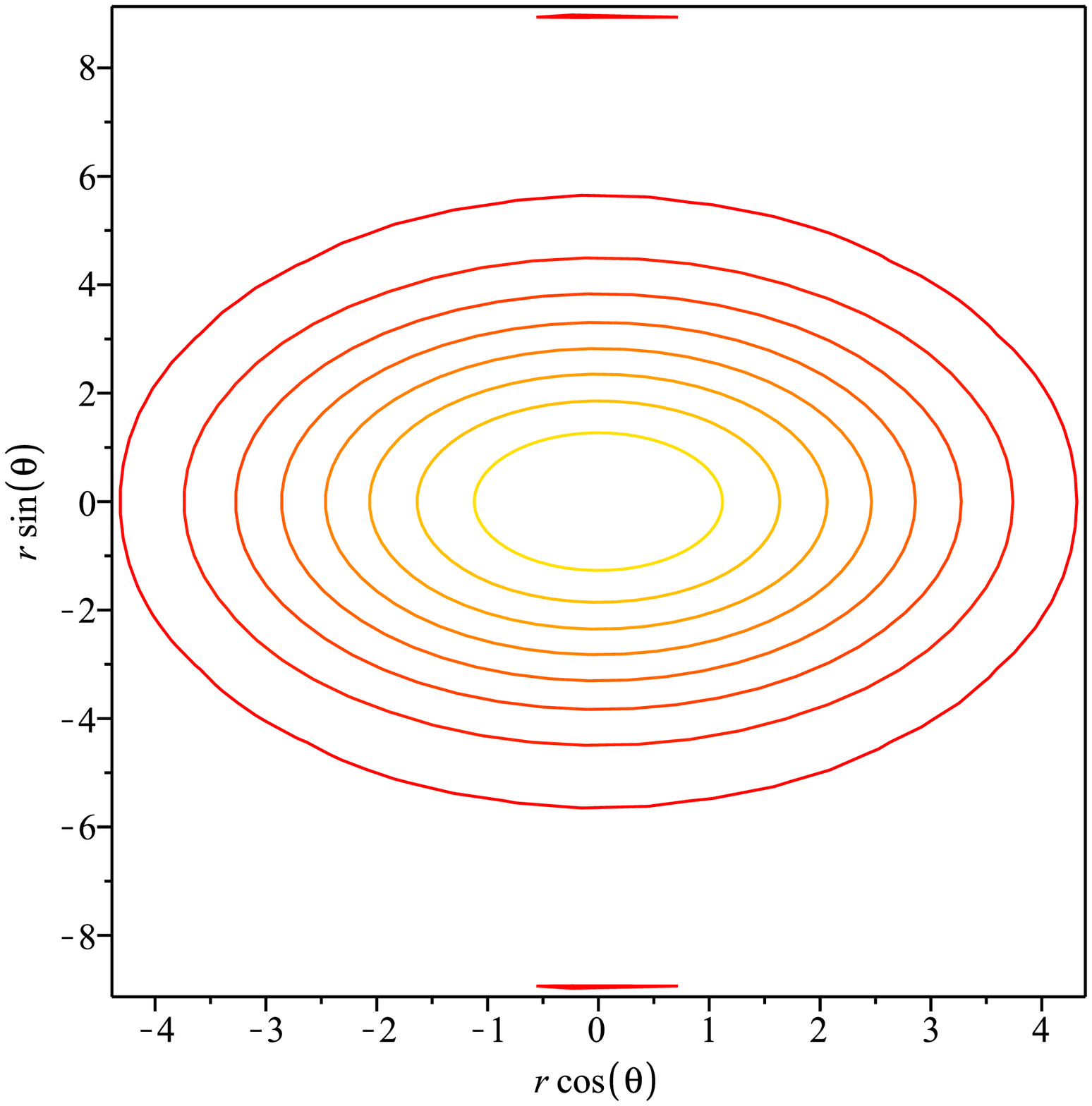,width=5.8cm}} \hspace{0.5cm}{\textbf{Fig.\,10a}\ \
 \small Contour plot of $B(r, \theta, \varphi)$ for\\${}^{}$\hspace{1.4cm}  the cylindrical coordinate system} \end{minipage}\quad\ \
\begin{minipage}{5.8cm}
{\psfig{figure=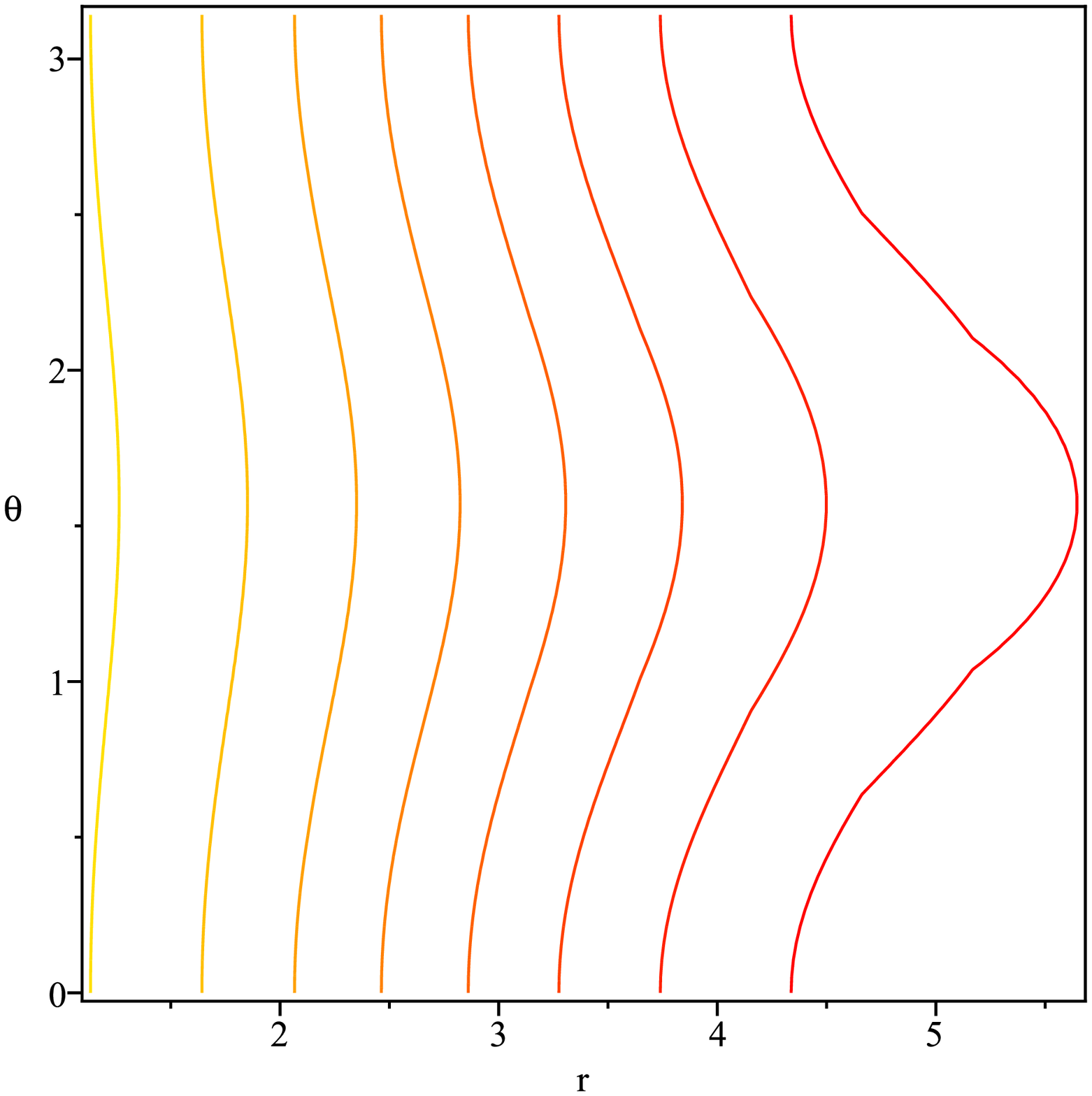,width=5.8cm}} \hspace{0.5cm}{\textbf{Fig.\,10b}\ \  \small Contour plot of $B(r, \theta, \varphi)$ in \\${}^{}$\hspace{1.9cm}$(r, \theta)$  coordinates} \end{minipage}
\begin{center}\begin{minipage}{5.8cm}
 {\psfig{figure=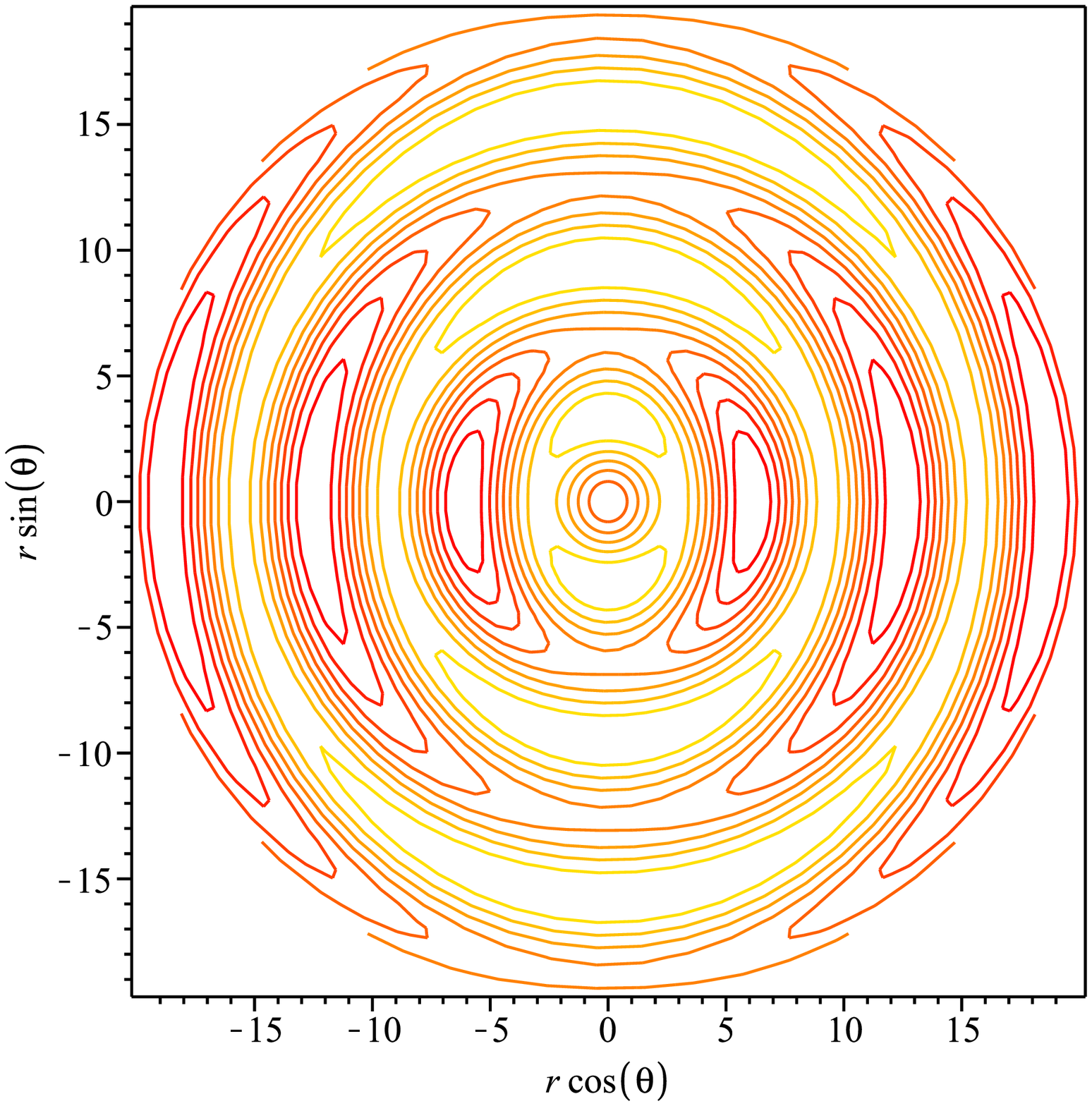,width=5.8cm}} \hspace{0.4cm}{\textbf{Fig.\,10c}\
 \small Contour plot of $r^2B(r, \theta, \varphi)$ for \\${}^{}$\hspace{1.4cm} the cylindrical coordinate system} \end{minipage}
 \end{center}
\end{example}
\section{Conclusions}

The paper studies  Abelian and Tauberian theorems for LRD random fields. We restrict our attention only to "classical" asymptotic results. Mercerian theorems, results for regularly varying functions with reminder, and other generalizations are beyond the scope of the paper. Some such extensions, formulated in terms of Hankel transforms, were developed in an excellent series of papers by \citealt{bin1,bin2,bin3,bin4} and can be used to investigate the fine asymptotic behavior of LRD random fields.

Using several new or less-known explicit examples of LRD covariance functions we illustrate and highlight the role of assumptions in the theorems to obtain Abelian and Tauberian results for different functional classes.

These examples can also be used in spatial statistics where one of the main problems is the development of covariance functions with some desirable properties (for example, LRD) in an explicit form  to fit experimental covariances. Some other models of random fields with singularities in the spectrum can be found in \citealt{hal, lav1} and \citealt{kly}.

\begin{acknowledgements}
The research of the first author was partially supported by the Commissions the European Communities grant  PIRSES-GA-2008-230804 (Marie Curie Action). The second author was partially supported by La Trobe University Research Grant "Stochastic Approximation in Finance and Signal Processing" and the Swedish Institute grant SI-01424/2007. The authors are grateful for the referee's comments, which helped to improve the style of the presentation.
\end{acknowledgements}

\end{document}